%% file: R-diag.tex
\documentclass[reqno,11pt]{amsart}
\usepackage{amsmath}
\usepackage{amsfonts,amssymb,amsmath,latexsym,wasysym,mathrsfs,bbm}
\usepackage[all]{xy}
\usepackage{color}
\usepackage{epsfig}

\def\R{\mathbb R}
\def\C{\mathbb C}

\def\N{\mathbb N}

\def\1{\mathbbm 1}

\def\e{\epsilon}

\def\tensor{\otimes}
\def\supp{\text{supp}\,}

\newcommand{\interior}[1]{\raise0.2ex\hbox{$\displaystyle{\mathop{#1}^{\circ}}$}}
\newcommand{\mx}[1]{\mathbf{#1}}

\renewcommand\phi{\varphi}

\newtheorem{theorem}{Theorem}[section]

\newtheorem{proposition}[theorem]{Proposition}
\newtheorem{definition}[theorem]{Definition}
\newtheorem{corollary}[theorem]{Corollary}

\newtheorem{lemma}[theorem]{Lemma}
\newtheorem{notation}[theorem]{Notation}
\numberwithin{equation}{section} 

\usepackage{palatino}

\def\A{\mathscr{A}}
\def\F{\mathbb F}

\def\H{\mathcal{H}}
\def\HH{\mathscr{H}}
\def\T{\vartheta}

\begin{document}

\title{Strong Haagerup inequalities for free $\mathscr{R}$-diagonal elements}
\author{Todd Kemp$^{\ast}$}
\address{$\ast$ Department of Mathematics \\ Cornell University \\ Ithaca, NY \\ 14850}
\email{tkemp@math.cornell.edu}
\author{Roland Speicher$^{\ast\ast}$}
\address{$\ast\ast$ Department of Mathematics and Statistics \\ Queen's University \\ Kingston, ON  \\  K7L 3N6}
\email{speicher@mast.queensu.ca}

\begin{abstract} In this paper, we generalize Haagerup's inequality \cite{Haagerup} (on convolution norm in the free group) to a very general context of $\mathscr{R}$-diagonal elements in a tracial von Neumann algebra; moreover, we show that in this ``holomorphic'' setting, the inequality is greatly improved from its originial form.  We give an elementary combinatorial proof of a very special case of our main result, and then generalize these techniques.  En route, we prove a number of moment and cumulant estimates for $\mathscr{R}$-diagonal elements that are of independent interest.  Finally, we use our strong Haagerup inequality to prove a strong ultracontractivity theorem, generalizing and improving the one in \cite{Biane 2}. \end{abstract}

\maketitle

\section{Introduction}

There is an interesting phenomenon which often occurs in holomorphic spaces.  A theorem in the context of a function space (for example a family of norm-estimates, such as the $L^p$-bound of the Riesz projection, \cite{Rudin}) takes on a stronger form when restricted to a holomorphic subspace.  For example, $L^p$-bounds often shrink, and have meaningful extensions to the regime $p<1$.  For our purposes, the most relevant example is {\em Janson's strong hypercontractivity theorem} \cite{Janson}, discussed below.  In algebraic terms, this theorem states that a certain semigroup has better properties when acting on the algebra generated by i.i.d.\  {\em complex} Gaussians than on the algebra generated by i.i.d.\ {\em real} Gaussians.  The latter is a $\ast$-algebra while the former is far from one; we will exploit this difference in what follows.

\medskip

In this paper, we will primarily be concerned with one prominent non-commutative norm inequality: the {\em Haagerup inequality}.  It first arose in \cite{Haagerup}, where it was the main estimate used to foster an example of a non-nuclear $C^\ast$-algebra with the metric approximation property.  In the context of that paper, Haagerup's inequality takes the following form:

\begin{theorem}[\cite{Haagerup}, Lemma 1.4]\label{classical Haagerup inequality} Let $\F_k$ be the free group on $k$ generators, and let $f\in \ell^2(\F_k)$ be a function supported on the subspace generated by words in $\F_k$ of length $n$.  Then $f$ acts as a convolutor on $\ell^2(\F_k)$, and its convolution norm $\|f\|_\ast = \sup_{\|g\|_2=1}\|f\ast g\|_2$ satisfies
\[ \|f\|_\ast \le (n+1)\|f\|_2. \]
\end{theorem}
Note that the convolution product is just the usual product in the von Neumann algebra generated by the left-regular representation of $\F_k$ (known as the {\em free group factor $L(\F_k)$}), and so in the language of operator algebras, the statement is that the (non-commutative) $L^2$-norm controls the operator norm on subspaces of uniform finite word-length, where the bound grows {\em linearly} with word-length.

\medskip

The Haagerup inequality, and its decendents, have played important roles in several different fields.  In the context of geometric group theory, the Haagerup inequality (and other constructions presented in \cite{Haagerup}) have evolved into {\em a-T-menability} or {\em property T} \cite{Valette 2}; in the context of Lie theory, Haagerup's inequality is related to {\em property RD} \cite{Lafforgue 2}.  It has proved useful for other operator algebraic applications: in \cite{Lafforgue 1}, Lafforgue uses the Haagerup inequality as a crucial tool in his proof of the Baum-Connes conjecture for cocompact lattices in $SL(3,\R)$; in this context, the precise order of growth of the Haagerup constant is immaterial (so long as it is polynomial).  On the other hand, the Haagerup inequality has proved useful in studying return probabilities and other statistics of random walks on groups (see \cite{CPS, Valette 1}), where the exact form of the Haagerup constant is important.

\medskip

Our main theorem, Theorem \ref{main theorem} below, is a strong Haagerup inequality in a general ``holomorphic'' setting -- i.e.\ a non-self-adjoint algebra.  In the special case of the free group factor, this amounts to considering  convolution operators which involve only generators of the group, not their inverses; the resulting Haagerup inequality (Corollary \ref{corollary free semigroup Haagerup inequality} below)  then has growth of order $\sqrt n$, where $n$ is the word-length.  

\medskip

There are two main approaches to norm estimates
in such a setting. A direct one (as used in the original approach
of Haagerup) is to work directly in the concrete representation
of the considered element as operator on a Hilbert space and
try to estimate the operator norm by considering the action
of the operator on vectors. A more indirect approach is
by recovering the operator norm as the limit of the $L^p$-norms as $p\to\infty$,
and therefore trying to get a combinatorial understanding of $L^p$-norms
for $p=2m$ even.  It is the latter approach which we take. 
Thus, we need a good (at least asymptotic) understanding of the moments of the involved
operators with respect to the underlying state. To our benefit, the moments of the generators
of free groups possess a lot of structure: namely the generators
are free in the sense of Voiculescu's free probability theory.

\medskip

Our strong Haagerup inequality is actually derived in a much more general setting: algebras generated by free $\mathscr{R}$-diagonal elements. We therefore handle not only the original framework of Haagerup (in the form of free Haar unitaries), but also free circular elements, and a wealth of other non-normal operators.

\medskip

There have been some predecessor of our strong Haagerup inequality
for the general $\mathscr{R}$-diagonal case. Namely, the one-dimensional
case was mainly addressed in \cite{Haagerup Larsen} and, in particular, in
\cite{Larsen}. Furthermore, \cite{Larsen} contains a very specialized multi-dimensional case, 
where the considered operator is a product of identically-distributed free $\mathscr{R}$-diagonal elements.
All these results relied on analytic techniques, using the theory of $\mathscr{R}$- and $\mathscr{S}$-transforms for probability measures on $\R$. However, in the genuine non-commutative case of polynomials in several non-commuting $\mathscr{R}$-diagonal elements, as we treat  it here, such analytical tools are unavailable to us, and so our analysis will rely
on the combinatorial machinery of free cumulants, as powered by free probability theory. 

\medskip

Our main tool is the moment-cumulant formula (Equation \ref{moments from cumulants}, below), which expresses the moments of the considered elements in a very precise combinatorial way in terms of free cumulants.  This allows
us to reduce the multi-dimensional case essentially to the one-dimensional case. (Note that this reduction is usually
the hardest part in such inequalities.) Whereas in some cases (as for circular elements) this reduction yields directly the desired result, in other cases -- namely when the cumulants of the $\mathscr{R}$-diagonal element may be negative (as it happens for
Haar unitaries, i.e., in the free group situation) -- we need an additional step.  Our strategy is to replace the original $\mathscr{R}$-diagonal element $a$ with a different $\mathscr{R}$-diagonal 
element $b$ whose cumulants are positive and dominate the
absolute values of the cumulants of $a$; this has to be
done in such a way that we have control over both the $L^2$-norm
and the operator norm of $b$ in terms of the corresponding
norms of $a$.  The technique we develop will, we hope, have more general applicability.

\medskip

Let us now give a precise definition of the arena for our Haagerup inequality.  Section \ref{section primer} contains brief introductions to all the terms used in what follows (and in the foregoing).

\begin{definition}\label{definition H(a,I)} Let $I$ be any indexing set, and let $\{a_i\,:\,i\in I\}$ be $\ast$-free identically distributed $\mathscr{R}$-diagonal elements in a $C^\ast$-probability space with state $\phi$; for convenience, let $a$ be a fixed $\mathscr{R}$-diagonal element with the same distribution.  Define $\H(a,I)$ to be the norm-closed $($non-$\ast)$ algebra generated by the $a_i$.  For each $n\ge 0$, define $\H^{(n)}(a,I)$ as the Hilbert subspace of $L^2(\H(a,I),\phi)$ of all elements of the form
\[ T=\sum_{|\mx{i}|=n} \lambda_{\mx{i}} a_{\mx{i}}, \]
where $\mx{i} = (i_1,\ldots,i_n)\in I^n$, $\lambda_\mx{i}\in\C$, and $a_{\mx i} = a_{i_1}\cdots a_{i_n}$.  We refer to $\H^{(n)}(a,I)$ as the {\bf $n$-particle space} $($relative to $a,I)$.
\end{definition}

The motivation for considering the algebra $\H(a,I)$ comes from the first author's paper \cite{Kemp}, and \cite{Biane 1}.  If $c$ is a circular element, then $L^2(\H(c,I),\phi)$ is a free analogue of the {\em Segal-Bargmann space} of \cite{Bargmann} -- i.e. the space
$\H L^2(\HH,\gamma)$ of holomorphic functions on a Hilbert space $\HH$ of dimension $|I|$, square-integrable with respect to a certain Gaussian measure $\gamma$.  The Segal-Bargmann space is the framework for the complex wave representation of quantum mechanics.  It played an important role in the constructive quantum field theory program in the mid- to late-twentieth century.

\medskip

There is a natural operator, the {\em Ornstein-Uhlenbeck operator} or {\em number operator} $N$ on $L^2(\HH,\gamma)$, which is related to the energy operator in quantum field theory.  In the classical (Gaussian) context, the Ornstein-Uhlenbeck semigroup $e^{-tN}$ satisfies a regularity property called {\em hypercontractivity}: for $1<p\le r<\infty$ the semigroup $e^{-tN}$ is a contraction from $L^p(\HH,\gamma)$ to $L^r(\HH,\gamma)$ for large enough time $t$.  When $e^{-tN}$ is restricted to the Segal-Bargmann space and its holomorphic $L^p$ generalizations, the time to contraction is shorter, as shown in \cite{Janson} and generalized in \cite{Gross}.  This {\em strong hypercontractivity} demonstrates that contraction properties of the Ornstein-Uhlenbeck semigroup improve in the holomorphic category.

\medskip

In \cite{Biane 2}, Biane showed how to canonically generalize the Ornstein-Uhlenbeck operator to the setting of free group factor, and proved that the resulting semigroup $e^{-tN_0}$ is hypercontractive.  He further showed that the semigroup $e^{-tN_0}$ satisfies an even stronger condition called {\em ultracontractivity}: it continuously maps $L^2$ into $L^\infty$ for {\em all} $t>0$, and for small time $\|e^{-tN_0}\|_{2\to\infty}$ is of order $t^{-3/2}$.  This result was proved using a version of the Haagerup inequality presented in \cite{Bozejko 1}.  We should note that, although this result is for the free group factor, the $n$-particle spaces used in the proof are {\em not} the same as in Theorem \ref{classical Haagerup inequality}, but are rather defined in terms of a generating family of {\em semicircular elements} defined in Section \ref{section primer}; nevertheless, the relevant Haagerup inequality {\em can} be proved from Theorem \ref{classical Haagerup inequality} using a central limit approach similar to the one in \cite{Voiculescu Dykema Nica}.

\medskip

It is Biane's free ultracontractivity theorem, along with our intuition that norm-inequalities improve in holomorphic categories, that motivated us to consider the same type of Haagerup inequality for $\mathscr{R}$-diagonal elements.  In the special case of circular elements, the first author showed in \cite{Kemp} that, as in the Gaussian case, in the holomorphic category -- in this case the spaces $L^p(\H(c,I),\phi)$ -- Biane's hypercontractivity result is trumped by Janson's strong hypercontractivity.  The first author further spelled out precisely the holomorphic structure inherent in $\H(c,I)$.  Our interpretation of $\mathscr{R}$-diagonal elements as ``holomorphic'' is more vague. Nevertheless, the algebra $\H(a,I)$ is a triangular algebra much like the space of bounded Hardy functions $H^\infty$ is (as a Banach algebra acting on $L^2(S^1)$).  More importantly, the kinds of norm estimates used in \cite{Kemp} have natural analogues for $\mathscr{R}$-diagonal elements.

\medskip

The following theorem, which is our strong version of Haagerup's inequality in the general $\mathscr{R}$-diagonal setting, is the main result of this paper.

\begin{theorem}\label{main theorem} Let $a$ be an $\mathscr{R}$-diagonal element in a $C^\ast$-probability space. There is a constant $C_a<\infty$ such that for all $T\in \H^{(n)}(a,I)$,
\begin{equation}\label{strong Haagerup} \|T\| \le C_a\sqrt{n}\, \|T\|_2. \end{equation}
In general, $C_a$ may be taken $\le 2^{10} \sqrt{e}\, \|a\|^2 / \|a\|_2^2$; if $a$ has non-negative free cumulants, $C_a$ may be taken $\le \sqrt{e}\, \|a\| / \|a\|_2$.
\end{theorem}

As a very special case (where the $a_i$ are free Haar unitaries), we deduce the following surprising strong version of the classical Haagerup inequality (Theorem \ref{classical Haagerup inequality}).

\begin{corollary}\label{corollary free semigroup Haagerup inequality} Let $k\ge 2$, let $\F_k$ be the free group on $k$ generators, and let $\F^+_k\subset\F_k$ be the free {\em semigroup} $($i.e.\ the set of all words in the generators, excluding their inverses$)$.  If $f\in \ell^2(\F^+_k)\subset\ell^2(\F_k)$ is supported on words of length $n$, then $f$ acts $($via the left-regular representation on the {\em full group} $\F_k)$ as a convolutor, with convolution norm
\[ \|f\|_\ast \le 2^{10}\sqrt{e}\,\sqrt{n}\, \|f\|_2. \]
\end{corollary}

\medskip

This paper is organized as follows.  In section \ref{section primer}, we give a brief introduction to free probability theory and $\mathscr{R}$-diagonal elements, in addition to setting the standard notation we will use throughout the paper.  In Section \ref{Circular Elements}, we provide a concrete bijection in order to calculate the moments of a circular element $c$; in it we derive, using more elementary techniques, a formula for $\|c^n\|$, confirming results in \cite{Oravecz} and \cite{Larsen}.  We then use this calculation, together with more involved combinatorial techniques, to estimate the norm of an element in the $n$-particle space $\H^{(n)}(c,I)$ for arbitrary indexing set $I$, and thus prove a special case of Theorem \ref{main theorem} in the circular context.

\medskip

In Section \ref{R-diagonal Elements}, we show how to modify the techniques in Section \ref{Circular Elements} to prove Theorem \ref{main theorem} in general.  In the process, we derive bounds on the growth of the free cumulants of $\mathscr{R}$-diagonal elements and, given an $\mathscr{R}$-diagonal $a$, show how to construct another $\mathscr{R}$-diagonal element $b$ with all positive cumulants dominating the cumulants of $a$.  We also show that the Haagerup inequality affiliated to the space $\H L^2(\nu_a)$ of holomorphic functions square integrable with respect to the Brown measure $\nu_a$ of $a$ is consistent with Theorem \ref{main theorem}, which shows that $\nu_a$ does carry some information about the mixed moments of $a$.  Finally, in Section \ref{Strong Ultracontractivity}, we introduce a natural analogue of the Ornstein-Uhlenbeck semigroup affiliated with $\H(a,I)$, and prove a strong ultracontractivity theorem for it.

\section{A Free Probability Primer}\label{section primer}

In this section we collect all the relevant results from free probability theory that will be used in what follows.  Our descriptions will be brief, as this material is quite standard and is explained in depth in the book \cite{Nica Speicher Book}.

\subsection{$C^\ast$-probability spaces}\label{subsection measures} Let $\A$ be a unital $C^\ast$ algebra, and let $\phi$ be a faithful state on $\A$ (i.e.\ for $a\in\A$, $\phi(a^\ast a)$ only vanishes when $a=0$).  The pair $(\A,\phi)$ is a {\em $C^\ast$-probability space}.  Elements of $\A$ are {\em non-commutative random variables} (which we will often refer to simply as {\em random variables}).  (Some authors prefer to reserve the term `random variable' for self-adjoint elements; in our context, all elements of $\A$ are treated equally.)  The motivating example is afforded by the commutative von Neumann algebra $L^\infty(\Omega, \mathcal{F}, P)$ of a probability space.  It comes equipped with the faithful state $\phi = \int_\Omega \cdot\; dP$; the random variables in this context are bounded random variables in the usual sense.

\medskip

In classical probability theory, any random variable $X$ has a probability distribution $\nu_X$ -- a measure on $\C$ which, among other things, determines the moments of $X$:
\[ \int_\Omega {X(\omega)}^n \overline{X(\omega)}{}^m\,dP(\omega) = \int_{\C} z^n\overline{z}^m\,d\nu_X(z,\bar{z}). \]
In the case of a real random variable $X$, $\nu_X$ is supported in $\R$ and we have $\int X^n\,dP = \int_{\R} t^n d\nu_X(t)$.  At least in the case of bounded random variables, these moment conditions uniquely determine the distribution, which is a compactly-supported probability measure.  The same holds true for normal elements in a $C^\ast$-probability space -- if $a$ is normal then there is a unique probability measure $\nu_a$ on $\C$ which satisfies
\begin{equation}\label{distribution} \phi\left(a^n (a^\ast)^m\right) = \int_{\C} z^n\overline{z}^m\,d\nu_a(z,\bar{z}), \end{equation}
and the measure $\nu_a$ is compactly supported.  Indeed, $\supp \nu_a$ is the spectrum of $a$, and the measure can be constructed using the spectral theorem: $\nu_a = \phi\circ E^a$ where $E^a$ is the spectral measure of $a$ in $\A$.

\medskip

If $a$ is not a normal element, then there is no measure satisfying Equation \ref{distribution}; more generally, given two elements in $\A$ that do not commute, there is no measure which represents their joint probability distribution (this is one way to state the Heisenberg uncertainty principle).  In the case where $(\A,\phi)$ is a tracial $W^\ast$-probability space ($\A$ is a von Neumann algebra, $\phi$ is a faithful normal tracial state) however, there is a best-approximation of a probability distribution called the {\em Brown measure}, introduced in \cite{Brown}.  If $a$ is normal, then its Brown measure coincides with its spectral measure, and so the Brown measure is also denoted $\nu_a$.  The Brown measure of $a$ always satisfies the moment condition $\phi(a^n) = \int_{\C} z^n\,d\nu_a(z,\bar{z})$, however it does not respect mixed-moments.

\subsection{The free group factors}\label{subsection free group factors} Free probability was invented by Voiculescu in \cite{Voiculescu} in order to import tools from classical probability theory into the study of {\em the free group factors} (specifically to address the still-open question of whether different free group factors are isomorphic).

\medskip

Let $k\ge 2$, and let $\F_k$ denote the free group on $k$ generators $u_1, u_2, \ldots, u_k$.  (We will also allow $k=\infty$ to denote the free group with countably-many generators.)  The {\em $k$th free group factor $L(\F_k)$} is the von Neumann algebra generated by the left-regular representation of $\F_k$ on $\ell^2(\F_k)$.  (Note:\ if $g\in\F_k$, then the image of $g$ in $L(\F_k)$ is an operator with $g^\ast = g^{-1}$.)  There is a natural state $\phi_k$ defined on $L(\F_k)$ induced by the function $g\mapsto \delta_{eg}$ on $\F_k$ (here $e$ is the identity in the group).  This state is faithful, normal, and tracial, making $(L(\F_k),\phi_k)$ into a $W^\ast$-probability space.

\medskip

There is a canonical representation of the free group factor on the full Fock space.  Let $\HH$ be a {\em real} Hilbert space, and let $\HH_\C = \C\tensor\HH$ be its complexification.  The {\em full Fock space} of $\HH$ is $\mathcal{F}(\HH) = \bigoplus_{j=0}^\infty (\HH_\C)^{\tensor j}$, where $\oplus$ and $\tensor$ are the Hilbert space direct sum and tensor product, and $(\HH_\C)^{\otimes 0}$ is defined to be the $\C$-span of an abstract vector $\Omega$ (not in $\H$) called the {\em vacuum vector}.

\medskip

For each $h\in\HH$, the {\em creation operator} $l(h)$ in $\mathscr{B}(\mathcal{F}(\HH))$ is uniquely defined by its action $l(h)(h_1\tensor\cdots\tensor h_j) = h\tensor h_1\tensor \cdots \tensor h_j$ on $(\HH_\C)^{\tensor j}$ (and $l(h)\Omega = h$).  The adjoint $l(h)^\ast$ is called the {\em annihilation operator}, and is given by $l(h)^\ast(h_1\tensor h_2\tensor \cdots \tensor h_j) = \langle h_1, h \rangle h_2\tensor \cdots \tensor h_j$ (and $l(h)^\ast\Omega=0$).  The operator $l(h)$ is not normal (if $h\ne 0$), but it is natural to consider the real part $X(h) = \frac{1}{2}(l(h)+l(h)^\ast)$.  For any $k$-dimensional real Hilbert space $\HH$, the von Neumann algebra generated by $\{X(h)\,:\,h\in\HH\}$ is isomorphic to $L(\F_k)$.  What's more, under this isomorphism, the state $\phi_k$ conjugates to the {\em vacuum expectation state} $\tau(X) = \langle X\Omega,\Omega \rangle$.

\medskip

Let $e_1,\ldots, e_k$ be an orthonormal basis for $\HH$.  The algebra $W^\ast\{X(h)\,:\,h\in\H\}\cong L(\F_k)$ is, of course, generated by the set $\{X(e_1),\ldots, X(e_k)\}$.  It is important to note that the isomorphism does {\em not} carry the generators $u_1,\ldots, u_k$ in $\F_k\subset L(\F_k)$ to the generators $X(e_1),\ldots,X(e_k)$.  Indeed, the two generating sets give two different, and important, families of non-commutative random variables: {\em Haar unitary} and {\em semicircular} elements, which we will discuss below.  In both cases, the relationship between different generators is a model of a non-commutative version of independence called {\em freeness}.

\subsection{Free cumulants and free independence} A normal random variable in a $C^\ast$-probability space is indistinguishable from a classical bounded complex random variable (indeed, one can construct a random variable with any given distribution $\nu$ as the identity function in the space $L^\infty(\nu)$.)  The important classical notion of independence of random variables, however, has no direct analog for pairs of non-commuting random variables.  The notion of {\em free independence} or {\em freeness}, introduced in \cite{Voiculescu} is a substitute, which is, in many ways, better.

\medskip

Let $\pi=\{V_1,\ldots,V_r\}$ be a partition of the set $\{1,\ldots, n\}$.  The partition is called {\em crossing} if for some $i\ne j$ there are numbers $p<q<p'<q'$ with $p,p'\in V_i$ and $q,q'\in V_j$. ({\em Notation}: we say $p\sim_\pi q$ if $p,q$ are in the same block of the partition $\pi$.  Thus, $\pi$ is crossing iff there are $p<q<p'<q'$ with $p\sim_\pi p'$, $q\sim_\pi q'$, and $p'\nsim_\pi q$.)  A {\em non-crossing partition} is one which is not crossing.  We represent a partition by connecting numbers in the same block $V_i$ of the partition.  The following figure gives four examples of non-crossing partitions of the set $\{1,\ldots, 6\}$.

\medskip

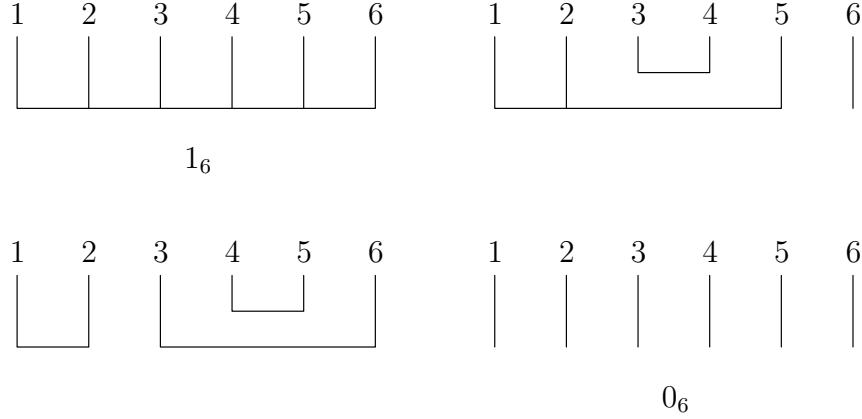
\begin{figure}[htbp]
\begin{center}
\input{fig1.pstex_t}
\caption{Four elements of $NC(6)$, including the minimal and maximal elements $0_6$ and $1_6$.}
\label{figure NC(6)}
\end{center}
\end{figure}

The set of non-crossing partitions of $\{1,\ldots, n\}$, denoted $NC(n)$, is partially-ordered under reverse refinement.  It is a lattice, in fact, with minimal element $0_n$ and maximal element $1_n$ as in Figure \ref{figure NC(6)}.  The M\"obius function $\mu_n$ of this lattice is well-known (see \cite{Kreweras}).  In particular, $\mu_n(0_n,1_n) = (-1)^{n-1}C_{n-1}$, where $C_n$ are the {\em Catalan numbers}
\begin{equation}\label{Catalan numbers} C_n = \frac{1}{n}\binom{2n}{n-1}. \end{equation}
More generally, for any $\sigma\in NC(n)$,
\begin{equation}\label{Moebius bound} |\mu_n(\sigma,1_n)|\le 4^{n-1}. \end{equation}
(The proof can be found contained in the proof of Proposition 13.15 in \cite{Nica Speicher Book}.)  It is worth noting that $C_n\le 4^n$ (and indeed $C^n\asymp 4^n$).

\medskip

Let $(\A,\phi)$ be a $C^\ast$-probability space.  Let $n>0$ and let $\pi$ be a partition in $NC(n)$.  For each block $V=\{i_1,\ldots, i_k\}$ in $\pi$, define the function $\phi_{V}\colon \A^n\to\C$ by $\phi_{V}[a_1,\ldots,a_n] = \phi(a_{i_1}\cdots a_{i_k})$.  Then define
$\phi_\pi\colon\A^n\to\C$ by $\phi_\pi[a_1,\ldots,a_n] = \prod_{V\in\pi} \phi_{V}[a_1,\ldots,a_n]$.  Finally, define the {\em free cumulants} of $(\A,\phi)$ to be the functionals $\{\kappa_\pi\,:\, \pi\in NC(n)\text{ for some }n>0\}$ by
\begin{equation}\label{cumulants} \kappa_\pi[a_1,\ldots,a_n] = \sum_{\stackrel{\scriptscriptstyle \sigma\in NC(n)}{\scriptscriptstyle \sigma\le\pi}} \phi_\sigma[a_1,\ldots, a_n]\,\mu_n(\sigma,\pi), \end{equation}
for each $\pi\in NC(n)$.  An immediate consequence of this definition is that the moments can be recovered from the free cumulants,
\[ \phi_\pi[a_1,\ldots,a_n] = \sum_{\stackrel{\scriptscriptstyle \sigma\in NC(n)}{\scriptscriptstyle \sigma\le\pi}} \kappa_\sigma[a_1,\ldots, a_n]. \]
(Indeed, this is the motivation for the inclusion of the coefficients $\mu_n(\sigma,\pi)$ in the definition of $\kappa_\pi$, for the M\"obius function is the convolution-inverse of the Zeta-function for the lattice $NC(n)$.)  As a special case, we have the formula \begin{equation}\label{moments from cumulants} \phi(a_1a_2\cdots a_n) = \sum_{\scriptscriptstyle \pi\in NC(n)} \kappa_\pi[a_1,\ldots,a_n]. \end{equation}

\medskip

Free cumulants allow a very easy statement of the definition of free independence, or freeness, of random variables.  Let $\kappa_n$ denote the free cumulant $\kappa_{1_n}$.  (These cumulants in fact contain all information about the cumulants, since all others can be built up block-wise by multiplication.)  Elements $a_1,\ldots, a_n$ in $\A$ are called {\em free} if, for $j\ge 2$ and $1\le i_1,\ldots,i_j\le n$, $\kappa_j[a_{i_1},\ldots,a_{i_j}]=0$ whenever there is at least one pair $1\le \ell,m\le j$ with $i_{\ell}\ne i_m$.  In other words, {\em random variables are free if all their mixed free cumulants vanish}.

\medskip

One can calculate that the generators $u_1,\ldots,u_n$ of $\F_n\subset L(\F_n)$ are free, as are the generators $X(e_1),\ldots, X(e_n)$ in the Fock-space representation of $L(\F_n)$; hence, this notion generalizes freeness from the free group context.  This approach mirrors the classical theory of cumulants in the method of moments (where the lattice considered is the lattice of {\em all} partitions).  All of the usual probabilistic constructions work: given any countable list of probability measures $\nu_j$, there is a $C^\ast$ probability space in which there are free random variables with distributions $\nu_j$ (one can construct the reduced free-product $C^\ast$ algebra of the $L^\infty(\nu_j)$, for example).

\subsection{$\mathscr{R}$-diagonal elements} As commented above, the operators $X(e_j)$ in the Fock-space representation of $L(\F_n)$ are semicircular elements: $s = X(e_j)$ has as distribution $\nu_s$ with
\[ d\nu_s(t) = \frac{1}{2\pi}\sqrt{4-t^2}\,dt. \]
Let $s_1,s_2$ be two free semicircular random variables.  The operator $c = (s_1+is_2)/\sqrt{2}$ (where $i = \sqrt{-1}$) is called a {\em circular element}.  It is non-normal, and so does not have a probability distribution.  (It's Brown measure is known, however, to be the uniform measure on the closed unit disc in $\C$.)  The $\ast$-cumulants  of a circular element (i.e.\ the free cumulants of tuples of operators all of the form $c$ or $c^\ast$) have a particularly nice form.  If $\varepsilon_j\in\{1,\ast\}$ then $\kappa_n[c^{\varepsilon_1},\ldots,c^{\varepsilon_n}]= 0$ for $n\ne 2$, and in fact only $\kappa_2[c,c^\ast] = \kappa_2[c^\ast,c] = 1$ are nonzero.

\medskip

Consider also a generator $u = u_j$ of $\F_k$.  Note that $\phi_k(u^n) = \delta_{n0}$, and the same holds true for $u^\ast = u^{-1}$.  The spectral measure of $u$ is thus the Haar measure on the unit circle, and such random variables are called {\em Haar unitary}.  The $\ast$-cumulants of a Haar unitary are not as restricted as those of a circular, but they follow a similar pattern.  The only nonvanishing cumulants $\kappa_n$ have $n$ even, and must have alternating $u$ and $u^\ast$ arguments:
\[ \kappa_{2n}[u,u^\ast,\ldots,u,u^\ast] = \kappa_{2n}[u^\ast,u,\ldots,u^\ast,u] = (-1)^{n-1}C_{n-1}, \]
the same as the M\"obius coefficents $\mu_n(0_n,1_n)$ of $NC(n)$ (and this is no coincidence).

\medskip

This connection between two widely known classes of non-selfadjoint random variables (circulars and Haar unitaries) motivated the second author, in \cite{Nica Speicher Fields Paper}, to introduce {\em $\mathscr{R}$-diagonal elements}.  A random variable $a$ in a $C^\ast$-probability space is $\mathscr{R}$-diagonal if its only novanishing cumulants are the alternating ones $\kappa_{2n}[a,a^\ast,\ldots,a,a^\ast]$ and $\kappa_{2n}[a^\ast,a,\ldots,a^\ast,a]$.  (The notation $\mathscr{R}$-diagonal derives from a characterization of such elements in terms of the multivariate {\em $\mathscr{R}$-transform}, a combinatorial free version of the logarithmic Fourier transform in classical probability theory.)

\medskip

Note that an $\mathscr{R}$-diagonal element's odd cumulants vanish.  (The term {\em even element} is used in this context, but is usually formulated in terms of {\em mixed} moments, so we do not use it for $\mathscr{R}$-diagonal elements.)  From Equations \ref{cumulants} and \ref{moments from cumulants} we see vanishing of odd cumulants is equivalent to vanishing of odd moments.  (A semicircular $s$ is even: its mixed moments are just its moments since it is self-adjoint, and, like a circular, only its second cumulant is nonzero: $\kappa_n[s,\ldots,s] = \delta_{n2}$.)  If $a$ is $\mathscr{R}$-diagonal, its {\em determining sequences} are $(\alpha_n[a])_{n=1}^\infty$ and $(\beta_n[a])_{n=1}^\infty$ defined by
\begin{equation}\label{determining sequences} \begin{aligned}
\alpha_n[a] &= \kappa_{2n}[a,a^\ast,\ldots,a,a^\ast], \\
\beta_n[a] &= \kappa_{2n}[a^\ast,a,\ldots,a^\ast,a]. \end{aligned}\end{equation}
If $a$ is in a tracial probability space (better yet if $\phi$ restricted to the algebra generated by $a$ and $a^\ast$ is tracial), then $\alpha_n[a] = \beta_n[a]$; in any case, these sequences contain all the information about the cumulants (and therefore mixed moments) of $a$ and $a^\ast$.

\medskip

$\mathscr{R}$-diagonal elements form a large class of (mostly) non-normal elements about which a great deal is known. In a sense, they are non-normal analogues of rotationally invariant distributions in $\mathbb{C}$; namely, the distribution of an $\mathscr{R}$-diagonal element is not changed if is multiplied by a free Haar unitary. This results in a special polar decomposition and relations with maximization problems for free entropy \cite{Nica Speicher Book,NSS,HP}. Our main theorem (\ref{main theorem}) supports the point of view that $\mathscr{R}$-diagonal elements can be considered as non-normal versions of holomorphic variables.
\medskip

Finally, we comment that there is a precise description of the Brown measure of an $\mathscr{R}$-diagonal element in terms of its $\mathscr{S}$-transform (another formal power-series associated to the moments of $a$).  The following theorem shows that $\mathscr{R}$-diagonal elements have rotationally-invariant Brown measures with nice densities.  Let $\times_p$ denote the polar Cartesian product (i.e.\ $[x,y]\times_p [0,2\pi)$ is the closed annulus with inner-radius $x$ and outer-radius $y$).
\begin{theorem}[Corollary 4.5 in \cite{Haagerup Larsen}]\label{theorem Brown measure of R-diag} If $a$ is $\mathscr{R}$-diagonal (and is not a scalar multiple of a Haar unitary), then its Brown measure $\nu_a$ is supported on
$\left(\|a^{-1}\|_2^{-1},\|a\|_2\right]\times_p\left[0,2\pi\right)$ if $a$ is invertible, and on the disc $\left[0,\|a\|_2\right]\times_p\left[0,2\pi\right)$ if it is not.  Moreover, $\nu_a$ is rotationally-invariant with density
\[ d\nu_a(r,\theta) = f(r)\,dr\,d\theta, \]
where $f$ is strictly positive on $\left(\|a^{-1}\|_2^{-1},\|a\|_2\right]$ or $\left[0,\|a\|_2\right]$ and has an analytic continuation to a neighbourhood of this interval in $\C$.
\end{theorem}

\section{Circular Elements}\label{Circular Elements}

In this section, we prove Theorem \ref{main theorem} in the special case that $a=c$ is circular.  Our proof in Section \ref{R-diagonal Elements} subsumes this one, but the techniques in this proof are new and interesting, and motivate the proof in what follows.  In Section \ref{c^n}, we give a new combinatorial proof that the $\ast$-moments of the powers of a circular element are the {\em Fuss-Catalan numbers}, defined in Equation \ref{Fuss-Catalan numbers} below.  (The main ideas of the construction in this section are due to Drew and Heather Armstrong, and we thank them for their contribution.)  In Section \ref{Haagerup for H(c,I)}, we use the asymptotics of the Fuss-Catalan numbers to demonstrate the strong Haagerup inequality for algebras generated by free circular elements.

\subsection{The powers of a circular element}\label{c^n} Let $c$ be a (variance $1$) circular element in a $C^\ast$-probability space $(\A,\phi)$.  The moments of $c^n$ were calculated first by Oravecz \cite{Oravecz} and Larsen \cite{Larsen}, each using a different approach to iterated free convolution of the $\mathscr{R}$-transform of $c$.  We will reproduce their results here, using more elementary combinatorial techniques.

\medskip

From Equation \ref{moments from cumulants}, we have
\begin{equation}\label{circular sum} \phi[(c^n(c^n)^\ast)^m] \hspace{0.1in} =  \sum_{\scriptscriptstyle \pi\in NC(2nm)}\hspace{-0.1in} \kappa_\pi[c_{n,m}], \end{equation}
where $c_{n,m}$ is the list
\begin{equation}\label{c list} c_{n,m} = \overbrace{\underbrace{c, \;\ldots , c,\,}_n\;\underbrace{c^\ast, \ldots, c^\ast,}_n \ldots, \underbrace{c, \;\ldots , c,\,}_n\;\underbrace{c^\ast, \ldots, c^\ast}_n}^{2m\text{ groups}}\;. \end{equation}
Since $c$ is circular, its only nonzero free cumulants are $\kappa_2[c,c^\ast]=1$ and $\kappa_2[c^\ast,c]=1$, hence the only nonzero terms in the above sum are those for which the partition $\pi$ is a pair partition $\pi\in NC_2(2mn)$ (each block is of size $2$), and for which each $c$ is paired to a $c^\ast$ in $c_{n,m}$.  We call such pairings {\em $\ast$-pairings}, and denote the set of $\ast$-pairings in $NC_2(2mn)$ by $NC_2^\ast(n,m)$.  Pictured below are two examples of elements in $NC_2^\ast(3,4)$.

\begin{figure}[htbp]
\begin{center}
\input{fig2.pstex_t}
\caption{Two $\ast$-pairings in $NC_2^\ast(3,4)$.}
\label{figure P*(3,4)}
\end{center}
\end{figure}
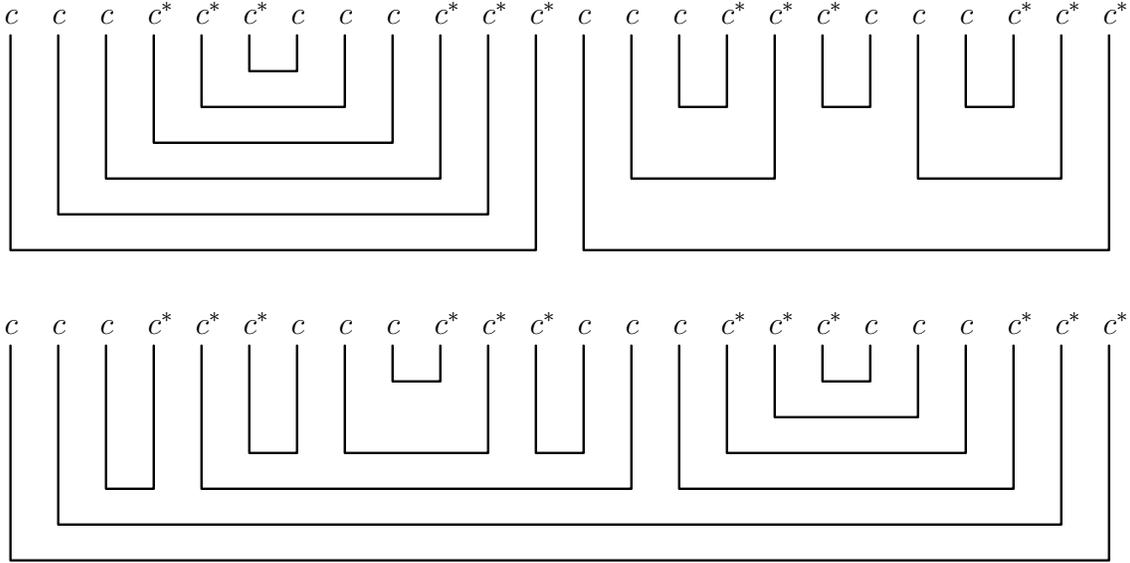

Since $\kappa_\pi[c_{n,m}] = 1$ whenever $\pi\in NC_2^\ast(n,m)$ and is $0$ otherwise, Equation \ref{circular sum} reduces to
\begin{equation}\label{m norm of c^n} \|c^n\|_{2m}^{2m} \hspace{0.1in} = \sum_{\pi\in NC_2^\ast(n,m)}\hspace{-.1in} 1 \hspace{0.1in} = \hspace{0.1in} |NC_2^\ast(n,m)|. \end{equation}
A non-crossing partition can be represented linearly as in Figures \ref{figure NC(6)} and \ref{figure P*(3,4)}, or equivalently on a circle, as seen below in Figure \ref{figure Round Table}.  As such, we can describe the problem of counting the elements in $NC_2^\ast(n,m)$ in the following medieval terms:

\medskip

{\em Knights and Ladies of the Round Table.}  King Arthur's Knights wish to bring their Ladies to a meeting of the Round Table.  There are $k = nm$ Knights (including Arthur himself) and each has one Lady.  Arthur wishes to seat everyone so that men and women alternate in groups of $n$, and in such a way that each Lady can converse with her Knight across the table without any conversations crossing.  How many possible seating plans are there?

\medskip

Letting $c$ stand for ``Knight'' and $c^\ast$ stand for ``Lady,'' the pictures in Figure \ref{figure Round Table} (which are the circular representations of the pairings from Figure \ref{figure P*(3,4)}) represent allowable seating plans.

\medskip

\begin{figure}[htbp]
\begin{center}
\input{fig3.pstex_t}
\caption{The $\ast$-pairings from Figure \ref{figure P*(3,4)}, in circular form.}
\label{figure Round Table}
\end{center}
\end{figure}
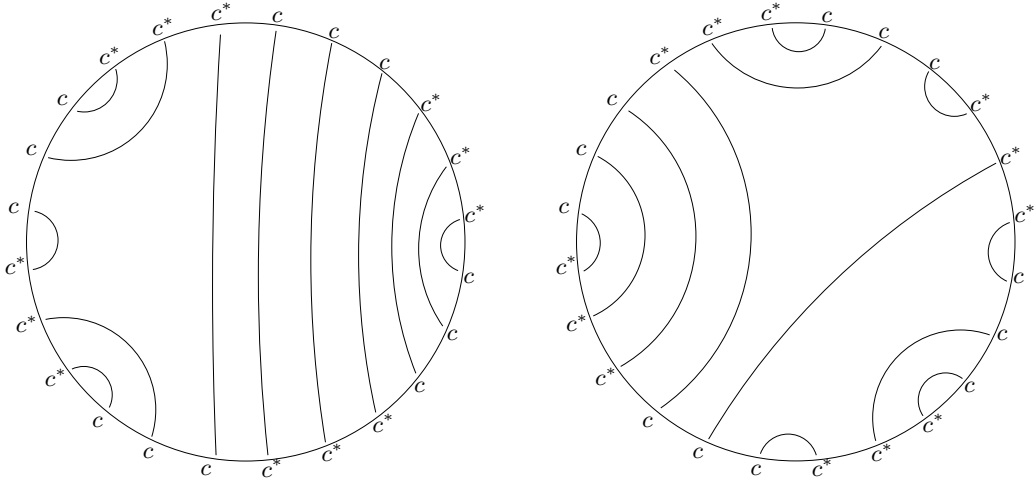

A related counting problem asks for pairings of the pattern $c_{n,m}$ where we relax the condition that each $c$ must be paired to a  $c^\ast$, but still required that {\em no two elements in a single $n$-block are paired together}.  Denote the set of all such non-crossing pairings as $\mathscr{T}(n,m)$ (so $NC_2^\ast(n,m)\subset\mathscr{T}(n,m)$).  As discussed in \cite{Biane Speicher}, this problem is the combinatorial counterpart to another moment problem, this time dealing with a semicircular element $s$.  Of course, since $s$ is selfadjoint, $(s^n(s^\ast)^n)^m = s^{2nm}$, and calculating these moments is routine.  Instead, the number of pairings in $\mathscr{T}(n,m)$ equals the moment $\phi(T_n(s)^{2m})$, where $T_n$ are the {\em Tchebyshev polynomials}.  While we do not have a nice schema for calculating $|\mathscr{T}(n,m)|$ explicitly (which we do for $|NC_2^\ast(n,m)|$ below), functional calculus for selfadjoint operators immediately yields that $\phi(T_n(s)^{2m})^{1/2m}\to n+1$ as $m\to\infty$ --- the norm $\|T_n(s)\|$ is {\em linear} in $n$, rather than in $\sqrt{n}$ as in Theorem \ref{main theorem} above.  This difference in size precisely reflects the improvement of Haagerup's inequality from $O(n)$ to $O(n^{1/2})$ behaviour for circular elements, and indeed for all $\mathscr{R}$-diagonal elements as discussed in Section \ref{R-diagonal Elements}.
\medskip

As to the {\em Knights and Ladies of the Round Table} problem, let us introduce some notation which will be useful throughout what follows.

\begin{notation}  Label the entries in $c_{n,m}$ with decreasing indices $n$ through $1$ in each block of $c$'s and increasing indices $1$ through $n$ in each block of $c^\ast$'s.
\begin{equation}\label{c c* numbering}  c_{n,m} = \underbrace{\mathop{c}_n, \; \mathop{c}_{n-1}, \; \ldots\, , \; \mathop{c}_2, \; \mathop{c}_1,} \; \underbrace{\mathop{c}_{1}{}^{\hspace{-0.02in}\ast}, \; \mathop{c}_{2}{}^{\hspace{-0.02in}\ast}, \; \ldots\, , \;  \mathop{c}_{n-1}{}^{\hspace{-0.1in}\ast}, \;  \mathop{c}_{n}{}^{\hspace{-0.03in}\ast}}, \; \ldots \, , \; \underbrace{\mathop{c}_n, \; \mathop{c}_{n-1}, \; \ldots \, , \; \mathop{c}_{2}, \; \mathop{c}_1}, \; \underbrace{\mathop{c}_1{}^{\hspace{-0.02in}\ast}, \;  \mathop{c}_{2}{}^{\hspace{-0.02in}\ast}, \; \ldots \, , \;  \mathop{c}_{n-1}{}^{\hspace{-0.1in}\ast}, \;  \mathop{c}_n{}^{\hspace{-0.03in}\ast}}
\end{equation}
We thus give each element of the list $c_{n,m}$ an address: $c(\ell,j)$ is the $\displaystyle{\mathop{c}_j}$ in the $\ell$th block of $c$'s, while $c^\ast(\ell,j)$ is the $\displaystyle{\mathop{c}_{j}{}^{\hspace{-0.02in}\ast}}$ in the $\ell$th block of $c^\ast$'s.
\end{notation}

\begin{lemma} For $1\le j \le n$, any $NC_2^\ast(n,m)$ must pair each $\displaystyle{\mathop{c}_j}$ to a $\displaystyle{\mathop{c}_{j}{}^{\hspace{-0.02in}\ast}}$.
\end{lemma}
\begin{proof}  The number of $c$'s between $c(\ell,j)$ and $c^\ast(\ell',j')$ is $n|\ell-\ell'|+j$, while the number of $c^\ast$'s between them is $n|\ell-\ell'|+j'$.  Let $\pi$ be a pairing which links $c(\ell,j)$ to (without loss of generality) $c(\ell',j')$ for some $j< j'$.  Since the number of $c$'s between $c(\ell,j)$ and $c^\ast(\ell',j')$ is greater than the number of $c^\ast$'s between them, $\pi$ must match at least one $c(k,i)$ between $c(\ell,j)$ and $c^\ast(\ell',j')$ to $c^\ast(k',i')$ where $k'<\min\{\ell,\ell'\}$ or $k'>\max\{\ell,\ell'\}$.  But then the blocks $\left\{c(\ell,j),c^\ast(\ell',j')\right\}$ and $\left\{c(k,i),c^\ast(k',i')\right\}$ in $\pi$ {\em cross}, and hence $\pi\notin NC_2(2nm)$.  Thus, $\pi\notin NC_2^\ast(n,m)$.
\end{proof}
\noindent We may note further that any non-crossing pairing which respects the labels in Equation \ref{c c* numbering} is, in fact, a $\ast$-pairing, and so enumerating $NC_2^\ast(n,m)$ amounts to counting the non-crossing pairings which respect those labels.  Using this observation, we proceed to define a bijection from $NC_2^\ast(n,m)$ to a set we can enumerate.

\begin{figure}[htbp]
\begin{center}
\input{fig4.pstex_t}
\caption{In the above $\ast$-pairing $\pi$, $1,3$ are $(\pi,3)$-connected, and $1,4$ are $(\pi^\ast,2)$-connected.}
\label{figure pi-connected}
\end{center}
\end{figure}
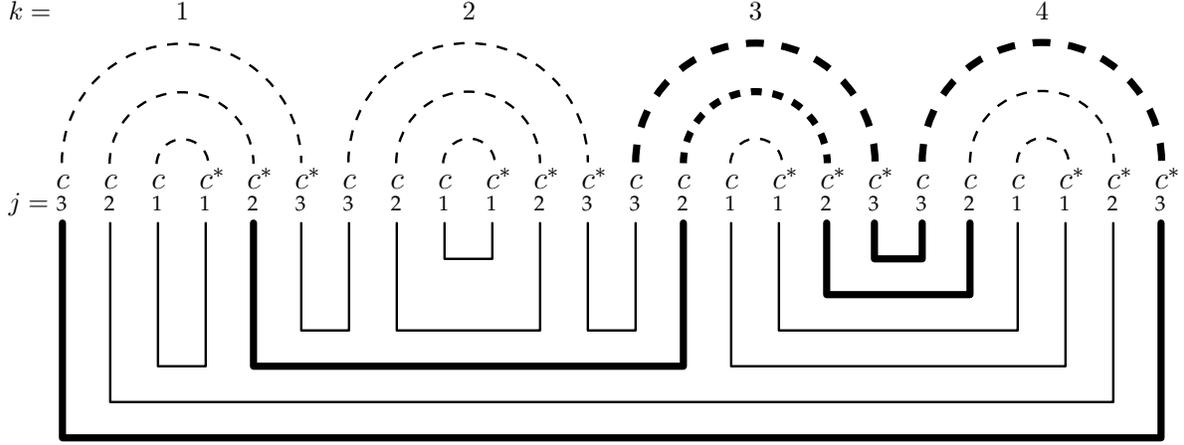

\begin{definition}\label{definition pi-connected}
Let $\pi\in NC_2^\ast(n,m)$, and let $1\le j \le n$. Say that $k,k'\in\{1,\ldots,m\}$ are {\bf $(\pi,j)$-connected} if there are $1\le k_1,\ldots,k_r\le m$ with $k_1>k$ such that $c(k,j)\sim_\pi c^\ast(k_1,j)$, $c(k_1,j)\sim_\pi c^\ast(k_2,j)$, \ldots, and $c(k_r,j)\sim_\pi c^\ast(k',j)$.  Similarly, say $k,k'$ are {\bf $(\pi^\ast,j)$-connected} if there are $1\le k_1,\ldots,k_r\le m$ with $k_1<k$ such that $c^\ast(k,j)\sim_\pi c(k_1,j)$, $c^\ast(k_1,j)\sim_\pi c(k_2,j)$, \ldots, and $c^\ast(k_r,j)\sim_\pi c(k',j)$.
\end{definition}
\noindent In other words, if we augment $\pi$ by connecting each pair $c(k,j),c^\ast(k,j)$, then $k,k'$ are $(\pi,j)$-connected if there is a(n initially increasing) path from $c(k,j)$ to $c^\ast(k',j)$ in the augmented pairing diagram; they are $(\pi^\ast,j)$-connected if there is a(n initially decreasing) path from $c^\ast(k,j)$ to $c(k',j)$.  If we exclude the conditions $k_1>k$ in $\pi$-connectedness and $k_1<k$ in $\pi^\ast$-connectedness, the two notions coincide (for example, $1$ and $4$ would be both $(\pi,2)$- and $(\pi^\ast,2)$-connected in Figure \ref{figure pi-connected}).  We find it convenient to treat them separately, however.

\medskip

We note the following, which is apparent in Figure \ref{figure pi-connected}.
\begin{lemma}\label{lemma monotone} If $k<k'$ are $(\pi,j)$-connected, then the sequence $k_1,\ldots,k_r,k'$ in definition \ref{definition pi-connected} is decreasing.  Likewise, if $k<k'$ are $(\pi^\ast,j)$-connected, then the sequence $k_1,\ldots,k_r,k'$ in definition \ref{definition pi-connected} is increasing.
\end{lemma}
\begin{proof} If $k<k'$ are $(\pi,j)$-connected, we have $k<k_1$, $c(k,j)\sim_\pi c^\ast(k_1,j)$ and $c(k_1,j)\sim_\pi c^\ast(k_2,j)$.  If $k_2 > k_1$, it follows that $c(k,j) < c(k_1,j) < c^\ast(k_1,j) < c^\ast(k_2,j)$, and hence there is a crossing.  The same argument applied at each pair $(k_\ell,k_{\ell+1})$ and at $(k_r,k')$ demonstrates the claim.  The argument for $(\pi^\ast,j)$-connectedness is similar.
\end{proof}

\begin{definition}\label{definition Phi} Given $\pi\in NC_2^\ast(n,m)$, define partitions $\Phi^\pi_n,\ldots,\Phi^\pi_1$ of $\{1,\ldots, m\}$ as follows: for $k,k'$ in $\{1,\ldots,m\}$, $k \sim_{\Phi^\pi_j} k'$ iff $k,k'$ are either $(\pi,j)$-connected or $(\pi^\ast,j)$-connected.
\end{definition}
\noindent That is, $\Phi^\pi_j$ is the image of $\left.\pi\right|_{\{c(k,j),c^\ast(j,k)\,:\,1\le k \le m\}}$ under the push-forward of the function $f_j$ from $\{c(j,k),c^\ast(j,k)\,;\,1\le k \le m\}$ to $\{1,\ldots,m\}$ which maps $c(j,k)$ and $c^\ast(j,k)$ to $k$.  (Note that $f_j$ is monotone.)

\begin{figure}[htbp]
\begin{center}
\input{fig5.pstex_t}
\caption{The partitions $\Phi_3,\Phi_2,\Phi_1$ corresponding to the two $\ast$-pairings in Figure \ref{figure P*(3,4)}.}
\label{figure Phi}
\end{center}
\end{figure}
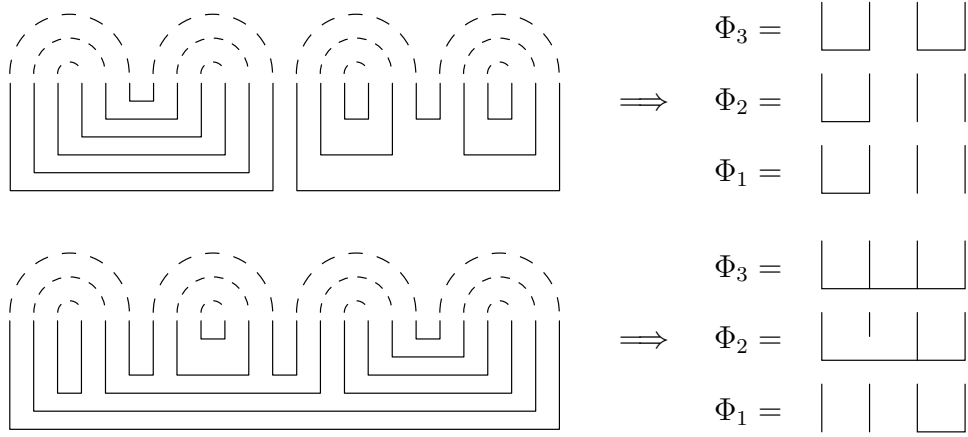

Figure \ref{figure Phi} shows the partitions $\Phi_j$ resulting from the $\ast$-pairings in \ref{figure P*(3,4)}; in it, we see that the $\Phi_j$ are non-crossing, and moreover they are refinement-decreasing -- in other words, they form a {\em multichain} (increasing sequence) in the lattice $NC(4)$: $\Phi_1\le\Phi_2\le\Phi_3$.  This holds generally for the $\Phi_j^\pi$ corresponding to any $\pi\in NC_2^\ast(n,m)$.

\begin{proposition}\label{proposition multichain} Let $\pi\in NC_2^\ast(n,m)$, and let $\Phi^\pi_1,\ldots,\Phi^\pi_n$ be the partitions in Definition \ref{definition Phi}.  Then the $\Phi^\pi_j$ are in $NC(m)$, and $\Phi^\pi_1\le\cdots\le\Phi^\pi_n$. \end{proposition}

\begin{proof} Since $f_j$ is monotone increasing and $\pi$ is non-crossing, $\Phi^\pi_j = (f_j)_\ast\left.\pi\right|_{\{c(k,j),c^\ast(j,k)\,:\,1\le k \le m\}}$ is non-crossing as well.  Now, let $1 < j \le n$, and suppose that $k<k'$ are connected by $\Phi^\pi_{j-1}$; thus, $k$ and $k'$ are either $(\pi,j-1)$-connected or $(\pi^\ast,j-1)$-connected.

\medskip

Suppose $k,k'$ are $(\pi,j-1)$-connected, and let $k_1,k_2,\ldots, k_r$ be a sequence connecting $c(k,j-1)$ to $c^\ast(k',j-1)$. By Lemma \ref{lemma monotone}, $k_1>k_2>\cdots>k_r>k'$.  Note that $c(k,j-1)<c^\ast(k,j)$, and so $c^\ast(k,j)$ must be paired to some $c(\ell_1,j)$ with $\ell_1>k$ -- otherwise $c(\ell_1,j) < c(k,j-1)< c^\ast(k,j) < c^\ast(k_1,j-1)$ resulting in a crossing.  If $\ell_1>k_1$ then there is a crossing at $c(k,j-1)<c^\ast(k,j)<c^\ast(k_1,j-1)<c(\ell_1,j)$; hence $\ell_1\le k_1$.  Suppose that $k' < \ell_1 < k_1$.  Then there is a $k_i$ with $k_{i+1}\le\ell_1<k_i$, giving a crossing with $c^\ast(k,j) < c^\ast(k_{i+1},j-1) < c(\ell_1,j) < c(k_i,j-1)$.  Hence, $k<\ell_1\le k'$.

\medskip

Inducting the previous argument, we find a chain $k< \ell_1 < \ell_2 < \cdots$ with $c^\ast(\ell_{i-1},j)\sim_\pi c(\ell_i,j)$, and each $\ell_i\le k'$.  Since there are only finitely many numbers between $k$ and $k'$, and since each $c^\ast(\ell_i,j)$ must be paired to a $c(\ell_{i+1},j)$ with $\ell_{i+1}>\ell_i$, it follows that $\ell_i = k'$ for some $i$.  Thus, $k,k'$ are $(\pi^\ast,j)$-connected.

\medskip

A similar argument shows that if $k<k'$ are $(\pi^\ast,j-1)$-connected then they are $(\pi,j)$-connected.  Hence, $\Phi^\pi_{j-1}$ is a refinement of $\Phi^\pi_j$, and so $\Phi^\pi_{j-1} \le \Phi^\pi_j$ in the lattice $NC(m)$.  \end{proof}

Denote by $NC^{(n)}(m)$ the set of all multichains of length $n$ in $NC(m)$.  Thus, Proposition \ref{proposition multichain} shows that the function $\mathscr{P}\colon\pi\mapsto(\Phi^\pi_1,\ldots,\Phi^\pi_n)$ is a map $NC_2^\ast(n,m)\to NC^{(n)}(m)$.  In what follows, we will show that $\mathscr{P}$ is a bijection.  To do so, we exhibit its inverse.

\medskip

To invert the above procedure for $\Phi\in NC(m)$, the idea (heuristically) is to ``fatten up'' each connecting line on the right-hand side of Figure \ref{figure Phi}, and assign pairings by ignoring the top connections (which identify each $c(k,j)$ with $c^\ast(k,j)$).

\begin{figure}[htbp]
\begin{center}
\input{fig6.pstex_t}
\caption{A ``fattened'' partition in $NC(8)$.}
\label{figure fattened}
\end{center}
\end{figure}
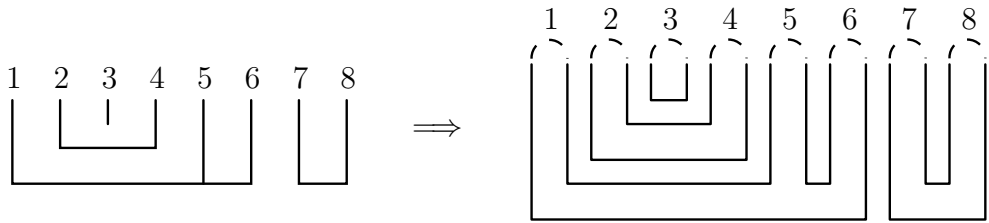

We can actually do this for each $j$ individually.

\begin{definition}\label{definition Phi inverse}Let $\Phi\in NC(m)$.  Define a partial-pairing $\pi^\Phi_j$ of $c_{n,m}$ as follows.  For each block $V = \{k_1<k_2<\cdots<k_r\}$ in $\Phi$, include in $\pi^\Phi_j$ the following pairings:
\[ c(k_1,j)\sim_{\pi^\Phi_j} c^\ast(k_r,j),\; c(k_r,j)\sim_{\pi^\Phi_j} c^\ast(k_{r-1},j),\; c(k_{r-1},j)\sim_{\pi^\Phi_j} c^\ast(k_{r-2},j),\; \ldots,\; c(k_2,j)\sim_{\pi^\Phi_j} c^\ast(k_1,j). \]
\end{definition}

\begin{proposition}\label{proposition P inverse well defined}
Given an $n$-multichain $\Phi_1\le\cdots\le\Phi_n$ in $NC^{(n)}(m)$, the pairing $\pi^{\Phi_1}_1\sqcup\cdots\sqcup \pi^{\Phi_n}_n$ is in $NC_2^\ast(n,m)$.
\end{proposition}
\noindent Note: the $\sqcup$'s above denote union of disjoint partial pairings.

\begin{proof} First, note that $\pi^{\Phi_j}_j$ is a refinement of the pull-back $f_j^\ast \Phi_j$, and so, again since $f_j$ is monotone and $\Phi_j$ is non-crossing, $\pi^{\Phi_j}_j$ is also non-crossing.  Let $k,k'$ be such that $c(k,j)\sim c(k',j)$ in $\pi^{\Phi_j}_j$, let $j'>j$ be such that $c(\ell,j')\sim c^\ast(\ell',j')$, and suppose there is a crossing between $\{c(k,j),c^\ast(k',j)\}$ and $\{c(\ell,j'),c^\ast(\ell',j')\}$.  There are eight possible arrangements -- we treat only the case $k < \ell < k' < \ell'$, and note the others may be treated similarly.  So, $c(k,j) < c(\ell,j') < c^\ast(k',j) < c^\ast(\ell',j')$.  Since $j'>j$, we have also $c(k,j) < c(\ell,j) < c^\ast(k',j) < c^\ast(\ell',j)$, and as $\Phi_{j}$ is a refinement of $\Phi_{j'}$, $\ell\sim_{\Phi_j}\ell'$ as well.  Thus there is a crossing in $\Phi_j$, which is a contradiction.  Hence, there are no crossings between $\pi^{\Phi_j}_j$ and $\pi^{\Phi_{j'}}_{j'}$ for any $1\le j<j' \le n$, and it follows that  $\pi^{\Phi_1}_1\sqcup\cdots\sqcup \pi^{\Phi_n}_n$ is in $NC_2(2nm)$.  By construction, it is a $\ast$-pairing, and so it is in $NC_2^\ast(n,m)$.  \end{proof}

Hence, the map $\mathscr{Q}\colon NC^{(n)}(m)\to NC_2^\ast(n,m)$ defined by $\mathscr{Q}(\Phi_1,\ldots,\Phi_n) =  \pi^{\Phi_1}_1\sqcup\cdots\sqcup \pi^{\Phi_n}_n$ is a well-defined function.  In fact, it is the inverse of $\mathscr{P}$.

\begin{proposition}\label{proposition bijection} The maps $\mathscr{P}\colon NC_2^\ast(n,m)\to NC^{(n)}(m)$ and $\mathscr{Q}\colon NC^{(n)}(m)\to NC_2^\ast(n,m)$ are inverses of each other.
\end{proposition}

\begin{proof} Let $\pi\in NC_2^\ast(n,m)$, and suppose that $c(k,j)\sim_{\pi}c^\ast(k',j)$.  Then $k,k'$ are in the same block of  $\Phi = \Phi^\pi_j$, and so by Definition \ref{definition Phi inverse}, $c(k,j)$ and $c^\ast(k',j)$ are connected in $\pi^{\Phi}_j$.  Hence, $\pi$ is a refinement of $\mathscr{Q}\circ\mathscr{T}(\pi)$.  On the other hand, suppose $c(k,j)$ and $c(k',j)$ are paired by $\mathscr{Q}\circ\mathscr{P}(\pi)$.  Then $\mathscr{P}(\pi) = (\Phi_1,\ldots,\Phi_n)$, where $k,k'$ are in the same block $V$ of $\Phi_j$, and moreover $k,k'$ are adjacent in the list $V = \{k_1,\ldots,k_r\}$ since, by Definition \ref{definition Phi inverse}, $\mathscr{Q}$ only creates pairings from adjacent elements of each block.  So, by definition \ref{definition Phi}, $k,k'$ are either $(\pi,j)$-connected or $(\pi^\ast,j)$-connected.  In either case, if the path connecting them were of length greater than $1$ then $k,k'$ would not be adjacent in the block $V$, since the sequence connecting them is monotone by Lemma \ref{lemma monotone}.  Hence, $k,k'$ are, in fact, connected in $\pi$.  This demonstrates that $\mathscr{Q}\circ\mathscr{P}(\pi)$ is a refinement of $\pi$, and so $\mathscr{Q}\circ\mathscr{P}(\pi) = \pi$.

\medskip

Now, let $(\Phi_1,\ldots,\Phi_n)\in NC^{(n)}(m)$, and let $(\Psi_1,\ldots,\Psi_n) = \mathscr{P}\circ\mathscr{Q}(\Phi_1,\ldots,\Phi_n)$.  If $k\sim_{\Phi_j} k'$ for $k<k'$, then there is a block $V$ of $\Phi_j$ including $k,k'$: $V = \{k_1 < \cdots < k_{r} < k < k_{r+1} < \cdots < k_s < k' < k_{s+1} < \cdots < k_t\}$.  Then $\pi^{\Phi_j}_j$ includes the pairings $c(k,j)\sim c^\ast(k_{r+1},j), \; c(k_{r+1},j)\sim c(k_{r+2},j), \; \ldots, \; c(k_s,j)\sim c^\ast(k',j)$; in particular, letting $\pi = \mathscr{Q}(\Phi_1,\ldots,\Phi_n)$, we have a path $(\pi,j)$-connecting $k$ and $k'$.  Hence, by Definition \ref{definition Phi}, $k\sim_{\Psi_j} k'$, and so $\Phi_j$ is a refinement of $\Psi_j$ for each $j$.  Conversely, if $k\sim_{\Psi_j} k'$, then $k,k'$ are either $(\pi,j)$-connected or $(\pi^\ast,j)$-connected.  Hence, there is a path connecting $k$ to $k'$ in $\pi^{\Phi_j}_j$, and so, by the action of $\mathscr{Q}$, $k$ and $k'$ must lie in the same block of $\Phi_j$ -- i.e.\ $k\sim_{\Phi_j} k'$.  This shows that $\Psi_j$ is a refinement of $\Phi_j$ for each $j$, and so we have shown that $\Psi_j = \Phi_j$ -- i.e.\ $\mathscr{P}\circ\mathscr{Q} = id_{NC^{(n)}(m)}$. \end{proof}

At this point, we have reproduced the results of Larsen using the above constructive approach.  The set $NC^{(n)}(m)$ is a well-studied combinatorial structure, and its enumeration was calculated by Edelman in \cite{Edelman}.  The next result follows.

\begin{corollary}\label{corollary Fuss-Catalan} For all positive integers $n$ and $m$, the number of $\ast$-pairings $|NC_2^\ast(n,m)|$ is equal to $|NC^{(n)}(m)| = C^{(n)}_m$, where $C^{(n)}_m$ are the {\em Fuss-Catalan numbers}
\begin{equation}\label{Fuss-Catalan numbers} C^{(n)}_m = \frac{1}{m}\binom{m(n+1)}{m-1}. \end{equation}
\end{corollary}
Note, in particular, that setting $n=1$ yields the Catalan numbers $C^{(1)}_m = C_m = \frac{1}{m}\binom{2m}{m-1}$ from Equation \ref{Catalan numbers}, which count the set $NC(m)$.  The Fuss-Catalan numbers were also computed in a similar context in \cite{Bisch Jones}, where the central objects of study, the {\em Fuss-Catalan algebras} (a generalization of the Temperly-Lieb algebras) are generated by diagrams like Figure \ref{figure Round Table}, and hence the dimensions of the algebras (the number of essentially different such diagrams) are the numbers $C^{(n)}_m$.

\subsection{The Haagerup inequality in $\H(c,I)$}\label{Haagerup for H(c,I)} From Equation \ref{m norm of c^n} and Corollary \ref{corollary Fuss-Catalan}, we have calculated the $2m$-norms of the powers of a circular element,
\begin{equation}\label{m norm of c^n 2}
\|c^n\|_{2m} = \left[C^{(n)}_m\right]^{1/2m} = \left[\frac{1}{m}\binom{m(n+1)}{m-1}\right]^{1/2m}. \end{equation}
In particular, the $2$-norm is $\|c^n\|_2 = 1$.  We can calculate the norm $\|c^n\|$ by taking the limit as $m\to\infty$, which may be computed using Stirling's formula.  The result is
\begin{equation}\label{norm squared} \|c^n\|^2 = \lim_{m\to\infty} \left[\frac{1}{m}\binom{m(n+1)}{m-1}\right]^{1/m} = \frac{(n+1)^{n+1}}{n^n} = \left(1+\frac{1}{n}\right)^n\,(n+1) \le e\,(n+1). \end{equation}

\medskip

Now, in line with Theorem \ref{main theorem}, consider the algebra $\H(c) = \H(c,\{1\})$, the norm-closed algebra generated by $c$.  In this case, the $n$-particle space $\H^{(n)}(c)$ is spanned by $c^n$, and hence Equation \ref{norm squared} immediately yields the following strong Haagerup inequality.
\begin{proposition}\label{H(c) Haagerup inequality} For $n\ge 0$ and $T\in\H^{(n)}(c)$,
\[ \|T\| \le \sqrt{e}\,\sqrt{n+1}\,\|T\|_2. \]
\end{proposition}

In fact, we can use similar techniques to achieve the same inequality for the algebra $\H(c,I)$ for any countable indexing set $I$.  This jump, from $1$ to many (even infinite) dimensions is usually the hardest part of such analyses; we will see below that the freeness does all the work for us. Note, the algebra $\H(c,I)$ is canonically isomorphic to the $0$-holomorphic space $\H_0(\HH_\C)$ in \cite{Kemp} and the free Segal-Bargmann space $\mathscr{C}_{hol}(\HH)$ in \cite{Biane 1}, where $\HH_\C$ is a complex Hilbert space of dimension $|I|$.

\medskip

Let $T\in\H^{n}(c,I)$, so that $T = \sum_{|\mx{i}|=n} \lambda_{\mx{i}} c_{\mx{i}}$ for some scalars $\lambda_{\mx{i}}\in\C$ satisfying a summability condition guaranteeing that $\|T\|_2<\infty$ (see Equation \ref{2 norm of circular T} below), where $c_{\mx{i}} = c_{i_1}\cdots c_{i_n}$.  By the definition of $\H(c,I)$, the generating elements $c_{i_k}$ are variance $1$ and $c_{i_k},c_{i_{k'}}$ are $\ast$-free whenever $i_k\ne i_{k'}$.  Then we have the following multinomial expansion for the $2m$th moment of $|T|$:
\begin{equation}\label{multinomial}\begin{aligned} \|T\|_{2m}^{2m} &= \; \phi[(T T^\ast)^m]  \\ &= \sum_{|\mx{i}(1)|=\cdots=|\mx{i}(m)|=n\atop|\mx{j}(1)|=\cdots=|\mx{j}(m)|=n} \lambda_{\mx{i}(1)}\cdots\lambda_{\mx{i}(m)}\overline{\lambda_{\mx{j}(1)}}\cdots\overline{\lambda_{\mx{j}(m)}}\,
\phi\left(c_{\mx{i}(1)}c_{\mx{j}(1)}^\ast\cdots c_{\mx{i}(m)}c_{\mx{j}(m)}^\ast\right). \end{aligned} \end{equation}
In particular, setting $m=1$,
\[ \|T\|_2^2 = \sum_{|\mx{i}|=|\mx{j}|=n} \lambda_{\mx{i}}\overline{\lambda_{\mx{j}}}\,\phi\left(c_{\mx{i}}c_{\mx{j}}^\ast\right). \]
The expression $\phi(c_{\mx{i}}c_{\mx{j}}^\ast)$ is a mixed moment of length $2n$, and can (by Equation \ref{moments from cumulants}) be expressed in terms of the cumulants of the $c_{\mx{i}}$:
\[ \phi\left(c_{\mx{i}}c_{\mx{j}}^\ast\right) = \sum_{\scriptscriptstyle \pi\in NC(2n)} \kappa_\pi[c_{i_1},\ldots,c_{i_n},c_{j_n}^\ast,\ldots,c_{j_1}^\ast]. \]
As the $c_i$ are circular (and so only the cumulants $\kappa_2[c,c^\ast] = \kappa_2[c^\ast,c] = 1$ are nonzero), only pair partitions $\pi$ which match $c$'s to $c^\ast$'s contribute to the sum.  Any such partition is in $NC_2^\ast(n,1)$, which contains only the partition $\varpi$

\begin{figure}[htbp]
\begin{center}
\input{nested.pstex_t}
\end{center}
\end{figure}

\noindent (the fact that there is only one follows from the calculation in Section \ref{c^n} that $|NC_2^\ast(n,1)| = C^{(n)}_1 = 1$).  So, we have
\begin{equation}\label{comma missing notation}  \|T\|_2^2 = \sum_{|\mx{i}|=|\mx{j}|=n} \lambda_{\mx{i}}\overline{\lambda_{\mx{j}}}\,\kappa_{\varpi}[c_{\mx{i}},c_{\mx{j}}^\ast]. \end{equation}
A note on notation: in Equation \ref{comma missing notation}, the $c_{\mx{i}}$ and $c^\ast_{\mx{j}}$ stand for {\em lists} of length $n$, not products of $n$ elements; i.e.\ there are implied commas.  We will use this convention whenever such expressions appear as arguments of cumulants in what follows.  To be clear, for the pairing $\varpi$ above, we have
\[ \kappa_{\varpi}[c_{\mx{i}},c_{\mx{j}}] = \kappa_{\varpi}[c_{i_1},\ldots,c_{i_n},c^\ast_{j_1},\ldots,c^\ast_{j_n}]
=\kappa_2[c_{i_1},c^\ast_{j_n}]\cdot\kappa_2[c_{i_2},c^\ast_{j_{n-1}}]\cdots\kappa_2[c_{i_n},c^\ast_{j_1}]. \]

\medskip

Now following Equation \ref{comma missing notation}, since the $c_{i_\ell}$ are $\ast$-free, $\kappa_\varpi[c_{\mx{i}},c_{\mx{j}}^\ast]=0$ unless each block of $\varpi$ contains like-indexed elements -- i.e.\ unless $\mx{i}=\mx{j}$, in which case $\kappa_\varpi = 1$.  Thus, we have the Pythagoreon formula
\begin{equation}\label{2 norm of circular T} \|T\|_2^2 = \sum_{|\mx{i}|=n} |\lambda_{\mx{i}}|^2. \end{equation}

\medskip

Following suit, for general $m>1$ we have
\[ \phi\left(c_{\mx{i}(1)}c_{\mx{j}(1)}^\ast\cdots c_{\mx{i}(m)}c_{\mx{j}(m)}^\ast\right)
= \sum_{\scriptscriptstyle \pi\in NC(2nm)} \kappa_\pi [c_{\mx{i}(1)},c_{\mx{j}(1)}^\ast,\ldots,c_{\mx{i}(m)},c_{\mx{j}(m)}^\ast]. \]
Once again, since the $c_{i(k)_\ell}$ are circular elements, the only partitions $\pi$ which contribute to the sum are those which pair $c$'s with $c^\ast$'s -- i.e.\ $\pi\in NC_2^\ast(n,m)$.  This, with Equation \ref{multinomial}, yields
\[ \|T\|_{2m}^{2m} = \sum_{\pi\in NC_2^\ast(n,m)} \sum_{|\mx{i}(1)|=\cdots=|\mx{i}(m)|=n\atop|\mx{j}(1)|=\cdots=|\mx{j}(m)|=n} \lambda_{\mx{i}(1)}\cdots\lambda_{\mx{i}(m)}\overline{\lambda_{\mx{j}(1)}}\cdots\overline{\lambda_{\mx{j}(m)}}\,
\kappa_\pi [c_{\mx{i}(1)},c_{\mx{j}(1)}^\ast,\ldots,c_{\mx{i}(m)},c_{\mx{j}(m)}^\ast]. \]
Many of the above terms are in fact $0$, since the $c_{i(k)_\ell}$ are $\ast$-free.  Indeed, the mixed cumulant $\kappa_\pi$ in the above sum is nonzero only when the indices of terms paired by $\pi$ are all equal (and in this case it is $1$).  We record this with the function $\delta(\pi,\mx{i}(1),\mx{j}(1),\ldots,\mx{i}(m),\mx{j}(m))$ defined to equal $0$ whenever $\pi$ pairs any $c_{i(k)_\ell}$ with a $c^\ast_{j(k')_{\ell'}}$ with $i(k)_\ell\ne j(k')_{\ell'}$, and $1$ if $\pi$ always pairs like-indexed $c$'s and $c^\ast$'s.  Thus
\[ \|T\|_{2m}^{2m} = \sum_{\pi\in NC^\ast(n,m)} \sum_{|\mx{i}(1)|=\cdots=|\mx{i}(m)|=n\atop|\mx{j}(1)|=\cdots=|\mx{j}(m)|=n} \lambda_{\mx{i}(1)}\cdots\lambda_{\mx{i}(m)}\overline{\lambda_{\mx{j}(1)}}\cdots\overline{\lambda_{\mx{j}(m)}}\,
\delta(\pi,\mx{i}(1),\mx{j}(1),\ldots,\mx{i}(m),\mx{j}(m)). \]
Now, let us re-index the above sum.  Denote the indices $\{i(1)_1,\ldots,i(m)_n\}$ by $p_1,\ldots,p_{nm}$, and let $\lambda(p_1,\ldots,p_{nm}) =  \lambda_{\mx{i}(1)}\cdots\lambda_{\mx{i}(m)}$.  Note, in any nonzero term in the above sum,
the indices appearing in the product $\overline{\lambda_{\mx{j}(1)}}\cdots\overline{\lambda_{\mx{j}(m)}}$ are exactly those paired to $p_1,\ldots,p_{nm}$ by $\pi$; identifying the pairing $\pi$ with its corresponding permutation, we then have
\begin{equation}\label{pre CS}
\|T\|_{2m}^{2m} = \sum_{\pi\in NC_2^\ast(n,m)}\sum_{p_1,\ldots,p_{nm}} \lambda(p_1,\ldots,p_{nm})\overline{\lambda(p_{\pi(1)},\ldots,p_{\pi(nm)})}. \end{equation}
Applying the Cauchy-Schwarz inequality to the interior summation yields, for each $\pi$,
\[ \begin{aligned} \sum_{p_1,\ldots,p_{nm}} \lambda(p_1,\ldots,p_{nm})&\overline{\lambda(p_{\pi(1)},\ldots,p_{\pi(nm)})} \\
\le &\left[\sum_{p_1,\ldots,p_{nm}} |\lambda(p_1,\ldots,p_{nm})|^2\right]^{1/2}\cdot\;
\left[\sum_{p_1,\ldots,p_{nm}} |\lambda(p_{\pi(1)},\ldots,p_{\pi(nm)})|^2\right]^{1/2}. \end{aligned} \]
Since the sum is over all $nm$-tuples of indices and $\pi$ is a permutation, the second term may be reordered to cancel the apparent $\pi$-dependence, yielding the same summation in both factors; i.e.\ the interior sum in Equation \ref{pre CS} is just
\[ \sum_{p_1,\ldots,p_{nm}} |\lambda(p_1,\ldots,p_{nm})|^2. \]
Returning to our original indexing scheme, this becomes
\[ \sum_{p_1,\ldots,p_{nm}} |\lambda(p_1,\ldots,p_{nm})|^2 = \sum_{|\mx{i}(1)|=\cdots=|\mx{i}(m)|=n} |\lambda_{\mx{i}(1)}\cdots\lambda_{\mx{i}(m)}|^2 = \sum_{|\mx{i}|=n}|\lambda_\mx{i}|^{2m} \le \left[\sum_{|\mx{i}|=n}|\lambda_\mx{i}|^2\right]^m, \]
and this last expression is $\|T\|_2^{2m}$ from Equation \ref{2 norm of circular T}.  Thus, Equation \ref{pre CS} and Corollary \ref{corollary Fuss-Catalan} together yield
\[ \|T\|_{2m}^{2m} \le \sum_{ \pi\in NC_2^\ast(n,m)}\|T\|_2^{2m} = C^{(n)}_m \|T\|_2^{2m}. \]
Taking $m$th roots and letting $m\to\infty$, referring to the same limit calculated in Equation \ref{norm squared}, we have thus proved the main theorem of this section:
\begin{theorem}\label{Theorem circular Haagerup} Let $c$ be a variance $1$ circular, and let $T\in\H^{(n)}(c,I)$ for some countable index set $I$.  Then
\[ \|T\| \le \sqrt{e}\,\sqrt{n+1}\,\|T\|_2. \]
\end{theorem}
We note that this inequality (with the $\sqrt{n+1}$ factor) bears some resemblance to what Bo\.zejko called {\em Nelson's inequality} in \cite{Bozejko 1}.  The context of his inequality is different, however (his estimate is for the creation and annihilation operators on the full Fock space separately), and our result cannot be derived from his.

\bigskip

\section{$\mathscr{R}$-diagonal Elements}\label{R-diagonal Elements}

In this section, we extend the techniques developed in Section \ref{Circular Elements} to all $\mathscr{R}$-diagonal elements.  A similar reduction of the multidimensional case to the one-dimensional case is possible, but there is an obstruction: the main argument goes through only when the mixed cumulants are non-negative.  We address this problem by replacing an $\mathscr{R}$-diagonal element with negative cumulants with a different $\mathscr{R}$-diagonal whose cumulants are positive and dominate the original's.

\medskip

In Section \ref{estimating moments}, we calculate the $2$-norm of an element $T$ in the $n$-particle space, and develop the main estimate (which generalizes the proof of Theorem \ref{Theorem circular Haagerup}) of higher moments of $|T|$ in terms of the absolute values of the cumulants.  Then, in Section \ref{Strong Haagerup inequalities}, we show how to replace a given $\mathscr{R}$-diagonal element with a different one who cumulants dominate the absolute values of the original's, and use this substitution to prove Theorem \ref{main theorem}.

\subsection{Estimating moments for $T\in\H^{(n)}(a,I)$}\label{estimating moments}

Let $a$ be an $\mathscr{R}$-diagonal element in a $C^\ast$-probability space, and let $T\in\H^{(n)}(a,I)$.  So, $T = \sum_{|\mx{i}|=n} \lambda_{\mx{i}} a_{\mx{i}}$ for some scalars $\lambda_{\mx{i}}\in\C$, where $\{a_i\,:\,i\in I\}$ are $\ast$-free $\mathscr{R}$-diagonal elements each with the same $\ast$-distribution as $a$.  As in Equation \ref{multinomial} above, we have the following multinomial expansion for the $2m$th moment of $|T|$:
\begin{equation}\label{a multinomial}\begin{aligned} \|T\|_{2m}^{2m} &= \; \phi[(T T^\ast)^m]  \\ &= \sum_{|\mx{i}(1)|=\cdots=|\mx{i}(m)|=n\atop|\mx{j}(1)|=\cdots=|\mx{j}(m)|=n} \lambda_{\mx{i}(1)}\cdots\lambda_{\mx{i}(m)}\overline{\lambda_{\mx{j}(1)}}\cdots\overline{\lambda_{\mx{j}(m)}}\,
\phi\left(a_{\mx{i}(1)}a_{\mx{j}(1)}^\ast\cdots a_{\mx{i}(m)}a_{\mx{j}(m)}^\ast\right). \end{aligned} \end{equation}
The term $\phi(a_{\mx{i}(1)}a_{\mx{j}(1)}^\ast\cdots a_{\mx{i}(m)}a_{\mx{j}(m)}^\ast)$ can be calculated, via Equation \ref{moments from cumulants}, as
\[ \phi\left(a_{\mx{i}(1)}a_{\mx{j}(1)}^\ast\cdots a_{\mx{i}(m)}a_{\mx{j}(m)}^\ast\right)
= \sum_{\pi\in NC(2mn)} \kappa_\pi [a_{\mx{i}(1)},a_{\mx{j}(1)}^\ast,\ldots,a_{\mx{i}(m)},a_{\mx{j}(m)}^\ast]. \]
Since the $a_{i(k)_\ell}$ are $\ast$-free, the above mixed cumulant is nonzero only when the indices of terms connected by $\pi$ are all equal.  We record this with the function $\delta(\pi,\mx{i}(1),\mx{j}(1),\ldots,\mx{i}(m),\mx{j}(m))$ defined above, which equals $0$ whenever $\pi$ connects two differently-indexed elements, and $1$ if all connected elements have like-indices.  It is, then, true that
\[ \phi\left(a_{\mx{i}(1)}a_{\mx{j}(1)}^\ast\cdots a_{\mx{i}(m)}a_{\mx{j}(m)}^\ast\right)
= \sum_{\pi\in NC(2mn)}\kappa_\pi
[a_{\mx{i}(1)},a_{\mx{j}(1)}^\ast,\ldots,a_{\mx{i}(m)},a_{\mx{j}(m)}^\ast]\delta(\pi,\mx{i}(1),\mx{j}(1),\ldots,\mx{i}(m),\mx{j}(m)).
\]
In the special case $m=1$, this reduces to
\begin{equation}\label{delta 2 norm sum} \phi\left(a_{\mx{i}}a_{\mx{j}}^\ast\right) = \sum_{\pi\in NC(2n)} \kappa_\pi [a_{\mx{i}},a_{\mx{j}}^\ast]\delta(\pi,\mx{i},\mx{j}). \end{equation}
Now, let $\pi$ be a partition with $\delta(\pi,\mx{i},\mx{j}) = 1$.  Thus, each block of $\pi$ connects only terms with a single index $i$.  Since $a_i$ is $\mathscr{R}$-diagonal, its only nonzero $\ast$-cumulants are $\kappa_{2n}[a_i,a_i^\ast,\ldots,a_i,a_i^\ast]$ and $\kappa_{2n}[a_i^\ast,a_i,\ldots,a_i^\ast,a_i]$.  Hence, $\pi$ still contributes a zero in Equation \ref{delta 2 norm sum} unless, in each block of $\pi$, the $a_i$'s and $a_i^\ast$'s alternate.  But in this case ($m=1$), all the $a_i^\ast$'s are to the right of all the $a_i$'s, and hence alternating sequences have length at most $2$.  So $\pi$ contributes only if it is a pair partition.  Since the cumulants $\kappa_2[a_j,a_j] = \kappa_2[a_j^\ast,a_j^\ast]=0$ for each $j$, such a $\pi$ only pairs $\ast$'s to non-$\ast$'s, and so $\pi$ is actually a $\ast$-pairing: $\pi\in NC_2^\ast(n,1)$.  As shown in Section \ref{Haagerup for H(c,I)}, the only element of $ NC_2^\ast(n,1)$ is $\varpi$.  So the sum in Equation \ref{delta 2 norm sum} reduces to at most a single term,
\[ \phi\left(a_{\mx{i}}a_{\mx{j}}^\ast\right) = \kappa_\varpi[a_{\mx{i}},a_{\mx{j}}^\ast]\delta(\varpi,\mx{i},\mx{j}). \]
Since $a_{\mx{i}}=a_{i_1}\cdots a_{i_n}$ and $a_{\mx{j}}^\ast = a_{j_n}^\ast\cdots a_{j_1}^\ast$, $\delta(\varpi,\mx{i},\mx{j})=1$ iff $\mx{i}=\mx{j}$, and in this case, $\kappa_\varpi[a_{\mx{i}},a_{\mx{i}}^\ast]$ is equal to the product $\kappa_2[a_{i_1},a^\ast_{i_1}]\cdots\kappa_2[a_{i_n},a^\ast_{i_n}]$ which (since the $a_i$ are identically distributed) equals $\kappa_2[a,a^\ast]^n$.  So Equation \ref{a multinomial} yields
\[ \|T\|_2^2 = \sum_{|\mx{i}|=|\mx{j}|=n} \lambda_{\mx{i}}\overline{\lambda_{\mx{j}}}\phi\left(a_{\mx{i}}a_{\mx{j}}^\ast\right)
= \sum_{|\mx{i}|=n} |\lambda_{\mx{i}}|^2 \kappa_2[a,a^\ast]^n. \]
Finally, we note that the second cumulant of a centred random variable is equal to its second moment (in general we may easily calculate that $\kappa_2[a,a^\ast] = Var(a)$), and since $\mathscr{R}$-diagonal elements have vanishing first moment, it follows that
\begin{equation}\label{2 norm of T a} \|T\|_2^2 = \sum_{|\mx{i}|=n} |\lambda_{\mx{i}}|^2 \|a\|_2^{2n}. \end{equation}

\medskip

Similar considerations are not enough to explicitly calculate higher moments, since alternating sequences can have greater length (e.g.\ in $\|T\|_4^4$, terms corresponding to partitions with blocks of sizes $2$ and $4$ may contribute), and calculations become unwieldy very quickly.  Nevertheless, we can estimate the higher norms using only pair partitions, to great effect.  In general, from Equation \ref{a multinomial} we have
\[ \|T\|_{2m}^{2m} = \sum_{\pi\in NC(2mn)} \sum_{|\mx{i}(1)|=\cdots=|\mx{i}(m)|=n\atop|\mx{j}(1)|=\cdots=|\mx{j}(m)|=n}  \lambda_{\mx{i}(1)}\cdots\lambda_{\mx{i}(m)}\overline{\lambda_{\mx{j}(1)}}\cdots\overline{\lambda_{\mx{j}(m)}}\cdot
\T[\pi,\mx{i}(1),\mx{j}(1),\ldots,\mx{i}(m),\mx{j}(m)], \]
where
\[ \T[\pi,\mx{i}(1),\mx{j}(1),\ldots,\mx{i}(m),\mx{j}(m)] = \kappa_\pi[a_{\mx{i}(1)},a_{\mx{j}(1)}^\ast,\ldots,a_{\mx{i}(m)},a_{\mx{j}(m)}^\ast]\delta(\pi,\mx{i}(1),\mx{j}(1),\ldots,\mx{i}(m),\mx{j}(m)). \]
Now, in any term where $\delta(\pi,\mx{i}(1),\mx{j}(1),\ldots,\mx{i}(m),\mx{j}(m))=1$, each block of $\pi$ connects only $a_i$'s and $a_i^\ast$'s for a single index $i$.  Since $a_i$ is $\mathscr{R}$-diagonal, its only nonvanishing $\ast$-cumulants are alternating, and so the term is zero unless $a$'s and $a^\ast$'s alternate within each block of $\pi$.  This is an important set of non-crossing partitions; we call it $NC^\ast(n,m)$ (so $NC_2^\ast(n,m)$ is the subset of $NC^\ast(n,m)$ consisting of only pair partitions).  It is important to note that, as per our definition of alternating, the size of each block of a partition in $NC^\ast(n,m)$ must be {\em even}.
(The sequence $a,a^\ast,\ldots,a,a^\ast,a$ is not alternating in our sense, since an $\mathscr{R}$-diagonal element still has vanishing cumulants for this list.)

\medskip

Using this notation, the above summation becomes
\[ \|T\|_{2m}^{2m} = \sum_{\pi\in NC^\ast(n,m)}\sum_{|\mx{i}(1)|=\cdots=|\mx{i}(m)|=n\atop|\mx{j}(1)|=\cdots=|\mx{j}(m)|=n} \lambda_{\mx{i}(1)}\cdots\lambda_{\mx{i}(m)}\overline{\lambda_{\mx{j}(1)}}\cdots\overline{\lambda_{\mx{j}(m)}}\cdot
\T[\pi,\mx{i}(1),\mx{j}(1),\ldots,\mx{i}(m),\mx{j}(m)]. \]
Fix $\mx{i}(1),\mx{j}(1),\ldots,\mx{i}(m),\mx{j}(m)$, and let $\pi\in NC^\ast(n,m)$ be such that $\delta(\pi,\mx{i}(1),\mx{j}(1),\ldots,\mx{i}(m),\mx{j}(m))=1$.  Let $\{V_1,\ldots,V_k\}$ be the blocks of $\pi$.  Since all indices of elements in a single block $V_j$ are equal (to, say, $i$), and since $a_i$ has the same distribution as $a$, we have that $\phi_{V_j}[a_{\mx{i}(1)},a_{\mx{j}(1)}^\ast,\ldots,a_{\mx{i}(m)},a_{\mx{j}(m)}^\ast] = \phi_{V_j}[a_{n,m}]$, where
\[ a_{n,m} = \overbrace{\underbrace{a, \;\ldots , a,\,}_n\;\underbrace{a^\ast, \ldots, a^\ast,}_n \ldots, \underbrace{a, \;\ldots , a,\,}_n\;\underbrace{a^\ast, \ldots, a^\ast}_n}^{2m\text{ groups}} \]
is independent of the indices.  Consequently, we have (for $\pi$ with $\delta(\pi,\mx{i}(1),\mx{j}(1),\ldots,\mx{i}(m),\mx{j}(m))=1$)
\begin{equation}\label{independent cumulants} \kappa_\pi[a_{\mx{i}(1)},a_{\mx{j}(1)}^\ast,\ldots,a_{\mx{i}(m)},a_{\mx{j}(m)}^\ast] = \kappa_\pi[a_{n,m}]. \end{equation}
Thus, for $\pi\in NC^\ast(n,m)$, we have
\[ \T[\pi,\mx{i}(1),\mx{j}(1),\ldots,\mx{i}(m),\mx{j}(m)] = \kappa_\pi[a_{n,m}]\delta(\pi,\mx{i}(1),\mx{j}(1),\ldots,\mx{i}(m),\mx{j}(m)), \] 
and so
\[ \|T\|_{2m}^{2m} = \sum_{\pi\in NC^\ast(n,m)} \kappa_\pi[a_{n,m}] \hspace{-0.12in} \sum_{|\mx{i}(1)|=\cdots=|\mx{i}(m)|=n\atop|\mx{j}(1)|=\cdots=|\mx{j}(m)|=n} \hspace{-0.2in} \lambda_{\mx{i}(1)}\cdots\lambda_{\mx{i}(m)}\overline{\lambda_{\mx{j}(1)}}\cdots\overline{\lambda_{\mx{j}(m)}}
\delta(\pi,\mx{i}(1),\mx{j}(1),\ldots,\mx{i}(m),\mx{j}(m)). \]
We now estimate this sum by associating to each $\pi\in NC^\ast(n,m)$ a refinement $\pi_r\in NC_2^\ast(n,m)$ as follows: for each block $V = \{k_1 < k_2 < \cdots < k_{2\ell}\}$ in $\pi$, the pairings $k_1\sim k_2$, $k_3\sim k_4$, \ldots, $k_{2\ell-1}\sim k_{2\ell}$ are in $\pi_r$.

\medskip

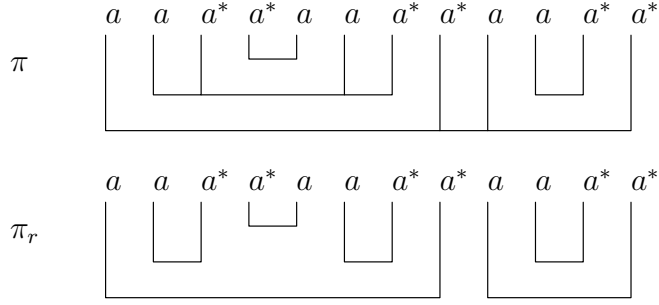
\begin{figure}[htbp]
\begin{center}
\input{fig7.pstex_t}
\caption{A partition $\pi\in NC^\ast(2,3)$, and the corresponding $\pi_r\in NC_2^\ast(2,3)$.}
\label{Figure splitting}
\end{center}
\end{figure}

Since $\pi_r$ is a refinement of $\pi$, if $\pi$ only connects like-indexed elements then $\pi_r$ does as well, and so $\delta(\pi,\mx{i}(1),\mx{j}(1),\ldots,\mx{i}(m),\mx{j}(m)) \le \delta(\pi_r,\mx{i}(1),\mx{j}(1),\ldots,\mx{i}(m),\mx{j}(m))$. Hence, we may estimate (by taking absolute values)
\[ \|T\|_{2m}^{2m} \le \sum_{\pi\in NC^\ast(n,m)} |\kappa_\pi[a_{n,m}]| \hspace{-0.2in} \sum_{|\mx{i}(1)|=\cdots=|\mx{i}(m)|=n\atop|\mx{j}(1)|=\cdots=|\mx{j}(m)|=n} \hspace{-0.2in} |\lambda_{\mx{i}(1)}\cdots\lambda_{\mx{i}(m)}\cdot\lambda_{\mx{j}(1)}\cdots\lambda_{\mx{j}(m)}|
\delta(\pi_r,\mx{i}(1),\mx{j}(1),\ldots,\mx{i}(m),\mx{j}(m)). \]
We can now reindex the interior sum the same way we did in Section \ref{Haagerup for H(c,I)}: denote the indices $\{i(1)_1,\ldots,i(m)_n\}$ by $p_1,\ldots,p_{nm}$, and this time let $\lambda(p_1,\ldots,p_{nm}) =  |\lambda_{\mx{i}(1)}\cdots\lambda_{\mx{i}(m)}|$.  Then allowing $\pi_r$ to refer both to the pair-partition and the associated permutation, we have
\[ \begin{aligned} &\sum_{|\mx{i}(1)|=\cdots=|\mx{i}(m)|=n\atop|\mx{j}(1)|=\cdots=|\mx{j}(m)|=n} \hspace{-0.2in} |\lambda_{\mx{i}(1)}\cdots\lambda_{\mx{i}(m)}\cdot\lambda_{\mx{j}(1)}\cdots\lambda_{\mx{j}(m)}|
\delta(\pi_r,\mx{i}(1),\mx{j}(1),\ldots,\mx{i}(m),\mx{j}(m)) \\
&=  \sum_{p_1,\ldots,p_{nm}} \lambda(p_1,\ldots,p_{nm})\lambda(p_{\pi_r(1)},\ldots,p_{\pi_r(nm)}) \\
&\le \left[\sum_{p_1,\ldots,p_{nm}}  \lambda(p_1,\ldots,p_{nm})^2\right]^{1/2}\cdot\left[ \sum_{p_1,\ldots,p_{nm}} \lambda(p_{\pi_r(1)},\ldots,p_{\pi_r(nm)})^2\right]^{1/2},
\end{aligned} \]
where we have applied the Cauchy-Schwarz inequality.  Since the sum is over all indices $p_1,\ldots,p_{nm}$ and since $\pi_r$ is a permutation, the second term above can be reindexed to yield the first term, and hence the interior sum is
\[ \le \sum_{p_1,\ldots,p_{nm}}  \lambda(p_1,\ldots,p_{nm})^2 = \sum_{|\mx{i}|=n} |\lambda_{\mx{i}}|^{2m} \le \left[\sum_{|\mx{i}|=n}|\lambda_{\mx{i}}|^2\right]^m. \]
Combining this with Equation \ref{2 norm of T a} yields the following estimate, which is the main lemma of this section.
\begin{lemma}\label{main lemma} Let $T \in \H^{(n)}(a,I)$ for $a$ $\mathscr{R}$-diagonal.  Then for $m\ge 1$,
\[ \|T\|_{2m} \le \left[\sum_{\pi\in NC^\ast(n,m)} |\kappa_\pi[a_{n,m}]|\right]^{1/2m} \frac{1}{\|a\|_2^n}\, \|T\|_2. \]
\end{lemma}
If the cumulants of $a$ are all non-negative, then $\kappa_\pi[a_{n,m}]\ge 0$ as well, and the above summation reduces to a one-dimensional calculation.

\begin{corollary}\label{corollary positive cumulants} If the cumulants of $a$ are non-negative, then
$\displaystyle{ \|T\| \le \frac{\|a^n\|}{\|a\|_2^n}\,\|T\|_2}$.
\end{corollary}

\begin{proof}  By Equation \ref{moments from cumulants},
\[ \|a^n\|_{2m}^{2m} = \phi([a^n(a^\ast)^n]^m) = \sum_{\pi\in NC(2nm)}\kappa_\pi[a_{n,m}]. \]
As explained above, since $a$ is $\mathscr{R}$-diagonal, $\kappa_\pi[a_{n,m}]=0$ unless $\pi\in NC^\ast(n,m)$. Thus, from Lemma \ref{main lemma}, we have
\[\begin{aligned} \|T\|_{2m}^{2m} &\le \left[\sum_{\pi\in NC^\ast(n,m)} |\kappa_\pi[a_{n,m}]|\right] \frac{1}{\|a\|_2^{2nm}}\, \|T\|_2^{2m} \\
&= \left[\sum_{\pi\in NC(2nm)} \kappa_\pi[a_{n,m}]\right] \frac{1}{\|a\|_2^{2nm}}\, \|T\|_2^{2m} = \frac{\|a^n\|_{2m}^{2m}}{\|a\|_2^{2nm}}\,\|T\|_2^{2m}. \end{aligned}\]
The result now follows by taking $2m$th roots, and letting $m$ tend to $\infty$. \end{proof}
Hence, in this case, the question of Haagerup's inequality is reduced to determining the growth-rate of $\|a^n\|/\|a\|_2^n$, which was addressed in \cite{Larsen} (and will be discussed in the next section).  However, if some cumulants of $a$ are negative, we must work harder to make such an estimate.

\subsection{Strong Haagerup inequalities}\label{Strong Haagerup inequalities}

To reduce the calculation in Section \ref{estimating moments} to the one-dimensional case when $a$ can have negative cumulants, our strategy is to replace $a$ with a different $\mathscr{R}$-diagonal element $b$ whose cumulants are positive and dominate the absolute values of $a$'s cumulants.  We will do this in a way that allows close control of both $\|b\|$ and $\|b\|_2$.  

\medskip

To begin, we bound the growth of the nonvanishing cumulants of $a$.

\begin{lemma}\label{lemma bound cumulants} Let $a$ be an $\mathscr{R}$-diagonal element in a $C^\ast$-probability space.  Then the nonvanishing cumulants of $a$ satisfy
\[ |\alpha_n[a]|,|\beta_n[a]| \le \frac{1}{2}(2^4\|a\|)^{2n}, \]
where $\alpha_n[a]$ and $\beta_n[a]$ are the determining sequences of $a$ from Equation \ref{determining sequences}. \end{lemma}
\begin{proof} From Equation \ref{cumulants}, we have
\[ \alpha_n[a] = \kappa_{2n}[a,a^\ast,\ldots,a,a^\ast] = \sum_{\sigma\in NC(2n)} \phi_\sigma[a,a^\ast,\ldots,a,a^\ast]\mu(\sigma,1_{2n}). \]
(The sum is over all of $NC(2n)$ since all $\sigma$ are less than $1_{2n}$, the largest element.)  Therefore, from Equation \ref{Moebius bound} we have
\[ |\alpha_n[a]| \le \sum_{\sigma\in NC(2n)} |\phi_\sigma[a,a^\ast,\ldots,a,a^\ast]|4^{2n-1}
= 4^{2n-1}\hspace{-.1in}\sum_{\sigma\in NC(2n)} \prod_{V\in\sigma} |\phi_{V}[a,a^\ast,\ldots,a,a^\ast]|. \]
Let $V_1,\ldots,V_r$ be the blocks of a given $\sigma\in NC(2n)$; so $|V_1|+\cdots+|V_r| = 2n$.  Well, $\phi_{V_j}[a,a^\ast,\ldots,a,a^\ast] = \phi(a^{\e_1}\cdots a^{\e_{|V_j|}})$ where $\e_i\in\{1,\ast\}$.  Since $\phi$ is a state on a $C^\ast$-algebra, this gives
\[ |\phi_{V_j}[a,a^\ast,\ldots,a,a^\ast]| \le \|a^{\e_1}\cdots a^{\e_{|V_j|}}\| \le \|a\|^{|V_j|}. \]
Hence, $|\phi_\sigma[a,a^\ast,\ldots,a,a^\ast]| = \prod_{j=1}^r |\phi_{V_j}[a,a^\ast,\ldots,a,a^\ast]| \le \prod_{j=1}^r \|a\|^{|V_j|} = \|a\|^{2n}$, and so
\[ |\alpha_n[a]| \le 4^{2n-1}\hspace{-.1in}\sum_{\sigma\in NC(2n)} \|a\|^{2n} = 4^{2n-1}C_{2n} \|a\|^{2n}. \]
The result for $\alpha_n[a]$ now follows from the fact that $C_{2n}\le 4^{2n}$.  The argument for $\beta_n[a]$ is identical. \end{proof}

Thus, we need only construct an $\mathscr{R}$-diagonal element whose determining sequences are positive and bounded below by $\frac{1}{2}(2^4\|a\|)^{2n}$.

\begin{lemma}\label{lemma bounding R-diag} Let $(\mathscr{A},\phi)$ be a $C^\ast$-probability space, and let $\gamma$ and $\lambda$ be positive constants. There exists an $\mathscr{R}$-diagonal element $b=b_{\gamma,\lambda}\in\mathscr{A}$ with $\alpha_n[b]=\beta_n[b]= \gamma\cdot\lambda^{2n}$. \end{lemma}

\begin{proof} As shown in \cite{Nica Speicher Book} (and also in \cite{Speicher Math. Anal. Paper}), there is a {\em free Poisson} element $p=p_{\frac{\gamma}{2},\lambda}$ which is self-adjoint and satisfies $\kappa_n[p,\ldots,p] = \frac{1}{2}\gamma\cdot\lambda^n$.  Let $p_1,p_2$ be free copies of this Poisson element, and let $q=p_1-p_2$.  As $\kappa_n$ is a linear combination of products of multilinear functionals $\phi_V$, and as $p_1$ and $-p_2$ are free (so their mixed cumulants vanish), we have
\[ \kappa_n[q,\ldots,q]  = \kappa_n[p_1,\ldots,p_1]+\kappa_n[-p_2,\ldots,-p_2]
= (1+(-1)^n)\kappa_n[p,\ldots,p] = \begin{cases} \gamma\cdot\lambda^n, & n\text{ even}, \\ 0, & n\text{ odd}\end{cases}. \]
Now, let $u$ be a Haar unitary $\ast$-free from $q$.  By Theorem 4.2(2) in \cite{Nica Speicher Duke Paper}, $b = qu$ is $\mathscr{R}$-diagonal.  (The conditions of the theorem require the $C^\ast$-probability space to be tracial; however, we may simply restrict $\phi$ to the unital $C^\ast$ algebra generated by the normal elements $q$ and $u$, where it is always a trace.)  Since $b$ is $\mathscr{R}$-diagonal, we can compute its determining sequences by
\[ \alpha_n[b]=\beta_n[b]=\sum_{\pi\in NC^\ast(n,1)}\phi_\pi[b,b^\ast,\ldots,b,b^\ast]\mu(\pi,1_{2n}). \]
Well, since $\pi\in NC^\ast(n,1)$, all blocks in $\pi$ are of even size and alternately connect $b$'s and $b^\ast$'s.  Hence, for each block $V$ in $\pi$,
\begin{equation}\label{b q cumulants} \phi_V[b,b^\ast,\ldots,b,b^\ast] = \phi[(bb^\ast)^{|V|/2}] = \phi[(quu^\ast q)^{|V|/2}] = \phi[q^{|V|}] = \phi_V[q,\ldots,q], \end{equation}
and thus $\phi_\pi[b,b^\ast,\ldots,b,b^\ast] =\phi_\pi[q,\ldots,q]$ for $\pi\in NC^\ast(n,1)$.

\medskip

Now, suppose $\sigma$ is a partition in $NC(2n)\setminus NC^\ast(n,1)$ -- i.e.\ $\sigma$ contains a block $V=\{k_1,\ldots,k_r\}$ with two successive elements $k_\ell<k_{\ell+1}$ of the same parity.  (Indeed, $NC^\ast(n,1)$ consists of non-crossing partitons whose blocks always successively pair $b$'s and $b^\ast$'s in the pattern $[b,b^\ast,\ldots,b,b^\ast]$ -- i.e.\ the blocks must alternately pair even and odd numbers in $\{1,\ldots,2n\}$.)  But then there is an odd number of elements between $k_\ell$ and $k_{\ell+1}$, and so some block in $\sigma$ must be of odd size.  Since $q$ is an even element, it follows that $\phi_\sigma[q,\ldots,q]=0$.  Hence, we also have
$\kappa_{2n}[q,\ldots,q] = \sum_{\pi\in NC^\ast(n,1)} \phi_\pi[q,\ldots,q]\mu_{2n}(\pi,1_{2n})$, and so from Equation \ref{b q cumulants},
\begin{align*} \alpha_n[b] &= \sum_{\pi\in NC^\ast(n,1)} \phi_\pi[b,b^\ast,\ldots,b,b^\ast]\mu_{2n}(\pi,1_{2n}) \\
&= \sum_{\pi\in NC^\ast(n,1)} \phi_\pi[q,\ldots,q]\mu_{2n}(\pi,1_{2n}) = \kappa_{2n}[q,\ldots,q] = \gamma\cdot\lambda^{2n}. \end{align*}  \end{proof}

\medskip

Following the argument of Corolloary \ref{corollary positive cumulants}, we see that if we choose an $\mathscr{R}$-diagonal element $b$ which satisfies $\alpha_n[b]\ge |\alpha_n[a]|$ for all $n$ then letting $b_{n,m}$ be the list corresponding to $[b^n (b^\ast)^n]^m$, we have $\kappa_\pi[b_{n,m}] \ge |\kappa_\pi[a_{n,m}]|$, and so
\[ \|b^n\|_{2m}^{2m} = \sum_{\pi\in NC^\ast(n,m)} \kappa_\pi[b_{n,m}] \ge \sum_{\pi\in NC^\ast(n,m)} |\kappa_\pi[a_{n,m}]|. \]
Hence, from Lemma \ref{main lemma}, we have
\begin{equation}\label{b over a estimate} \|T\|_{2m}\le \frac{\|b^n\|_{2m}}{\|a\|_2^n}\|T\|_2. \end{equation}
In order for this to yield useful information, we must choose $b$ in such a way that its variance and norm are well-controlled by those of $a$.  In the following lemma, we choose $b=b_{\gamma,\lambda}$ as in Lemma \ref{lemma bounding R-diag} to optimally bound the ratio $\|b^n\|/\|a\|_2^n$.

\begin{lemma}\label{lemma choose coefficients}  Let $a$ be $\mathscr{R}$-diagonal, and define $\lambda = 2^8\|a\|^2/\|a\|_2$ and $\gamma = \|a\|_2^2\lambda^{-2}$.  Set $b = b_{\gamma,\lambda}$, as in Lemma \ref{lemma bounding R-diag}.  Then
$\|b\|_2=\|a\|_2$, and
\[ \frac{\|b^n\|}{\|a\|_2^n} \le 2^{10}\,\sqrt{e}\,\sqrt{n}\, \frac{\|a\|^2}{\|a\|^2_2}. \]
\end{lemma}

\begin{proof} For $\mathscr{R}$-diagonal $b$, Corollary 3.2 in \cite{Larsen} says that $\|b^n\| \le \sqrt{e}\,\sqrt{n}\,\|b\|\,\|b\|_2^{n-1}$.  Note that, since $b$ is centred, $\|b\|_2^2 = \kappa_2[b,b^\ast]$ which, from Lemma \ref{lemma bounding R-diag}, equals $\gamma\cdot\lambda^2 = \|a\|_2^2$.  Hence,
\begin{equation}\label{direct ratio estimate} \frac{\|b^n\|}{\|a\|_2^n} = \frac{\|b^n\|}{\|b\|_2^n} \le \sqrt{e}\,\sqrt{n}\,\frac{\|b\|}{\|b\|_2} = \sqrt{e}\,\sqrt{n}\,\frac{\|b\|}{\|a\|_2}. \end{equation}
For the norm $\|b\|$, we have $b = qu$ where $u$ is unitary, and so $\|b\| = \|q\| = \|p_1-p_2\| \le 2\|p_1\|$.  The norm of a free Poisson was calculated in \cite{Voiculescu Dykema Nica}; the result is $\|p_1\| = \lambda(1+\sqrt{\gamma/2})^2$, so
\[ \sqrt{\gamma/2} = 2^{-1/2}\cdot\|a\|_2\lambda^{-1} = 2^{-8.5}\frac{\|a\|_2^2}{\|a\|^2} < 2^{-8.5}, \]
and so
\[ \|b\| \le 2\cdot2^8\frac{\|a\|^2}{\|a\|_2}\cdot(1+2^{-8.5})^2 \le 2^{10}\frac{\|a\|^2}{\|a\|_2}, \]
yielding the result. \end{proof}

We now stand ready to prove the main result of this paper.

\begin{proof}[Proof of Theorem \ref{main theorem}]
We will check that the element $b=b_{\gamma,\lambda}$ with coefficients chosen as in Lemma \ref{lemma choose coefficients} has all positive cumulants which dominate the absolute values of the cumulants of $a$.  First, we have (as used above) $\alpha_1[b]= |\alpha_1[a]|$.  For higher cumulants, using Lemma \ref{lemma bounding R-diag},
\[ \alpha_n[b] = \gamma\cdot\lambda^{2n} = \|a\|_2^2 \left(2^8\frac{\|a\|^2}{\|a\|_2}\right)^{2n-2} = \frac{1}{2}(2^4\|a\|)^{2n}\cdot
\left(\frac{\|a\|}{\|a\|_2}\right)^{2n-4}2^{8n-15}, \]
and since $n\ge 2$ and $\|a\|_2\le\|a\|$, this is $\ge \frac{1}{2}(2^4\|a\|)^{2n}$ which is, by Lemma \ref{lemma bound cumulants}, $\ge |\alpha_n[a]|$.  Having shown that $\alpha_n[b]\ge |\alpha_n[a]|$ for all $n$, we may now use  Equation \ref{b over a estimate}. We have (taking the limit as $m\to\infty$)
\[ \|T\| \le \frac{\|b^n\|}{\|a\|_2^n}\|T\|_2, \]
and from Lemma \ref{lemma choose coefficients} this yields the result:
\[ \|T\| \le 2^{10}\,\sqrt{e}\,\frac{\|a\|^2}{\|a\|^2_2}\,\sqrt{n}\,\|T\|_2. \]
If the cumulants of $a$ are all non-negative, then Equation \ref{b over a estimate} holds with $b=a$, and then Equation \ref{direct ratio estimate} yields the tighter estimate.
\end{proof}

Corollary \ref{corollary free semigroup Haagerup inequality} follows directly from Theorem \ref{main theorem}.  To be precise: if $u_1,\ldots,u_k$ are generators of $\F_k$, then the inclusions of $u_1,\ldots,u_k$ into $L(\F_k)$ are free Haar unitaries in the free group factor $L(\F_k)$ (this is discussed in Section \ref{subsection free group factors}).  The set of functions $g\in \ell^2(\F_k)$ supported on words in the $u_j$ (excluding their inverses) of length $n$ is equal to the $n$-particle space $\H^{(n)}(u,I_k)$ ($I_k=\{1,\ldots,k\}$) in the $W^\ast$-probability space $(L(\F_k),\phi_k)$, and a short calculation verfies that the norm on $\ell^2(\F_k)$ equals the norm in $L^2(L(\F_k),\phi_k)$.  Finally, the convolution norm is defined by $\|g\|_\ast = \supp_{f\ne 0} \|g\ast f\|_2/\|f\|_2$, which is the definition of the norm in the von Neumann algebra $L(\F_k)$.  So, Corollary \ref{corollary free semigroup Haagerup inequality} is indeed a special case of Theorem \ref{main theorem}.

\medskip

Note, the proof of Lemma \ref{lemma choose coefficients} actually produces a constant involving $2^9(1+2^{-8})^2 = 516.0078125$, far less than the stated $2^{10}=1024$.  However, since it is highly doubtful that this constant is optimal, there is little point quibbling.  That there is a constant at all -- i.e.\ that the behaviour is $O(n^{1/2})$ rather than $O(n)$, is the important, and surprising, fact.

\medskip

We also note that the sharp constant for $a$ with negative cumulants is {\em greater} than the sharp constant $\sqrt{e}\,(\|a\|/\|a\|_2)\,\sqrt{n}$ which holds when $\alpha_n[a]\ge 0$.  For example, consider a Haar unitary $u$, and the corresponding algebra $\H(u,\N)$.  For $k>1$ in $\N$, the element $T_k=u_1+\cdots+u_k$ is in the $1$-particle space, and satisfies $\|T_k\|_2=\sqrt{k}$ (Equation \ref{2 norm of T a}) and $\|T_k\|=2\sqrt{k-1}$ (as calculated in \cite{Haagerup Larsen}).  Thus
\[ \frac{\|T_k\|}{\|T_k\|_2} = 2\cdot\sqrt{\frac{k-1}{k}}. \]
Thus, if the Haagerup inequality $\|T\|\le C\, \|T\|_2$ (note $\|u\| = \|u\|_2=1$) holds for all $T\in\H^{(1)}(u,\N)$, then $C\ge 2 > \sqrt{e}$. It may be that $2\sqrt{n}$ is the optimal constant for $\H(u,\N)$, but we are as yet unable to calculate norms of elements in these $n$-particle spaces for $n>1$.

\medskip

We conclude this section with a discussion of Brown measure.

\begin{theorem}\label{theorem Brown measure Haagerup inequality} Let $a$ be an $\mathscr{R}$-diagonal element which is not a scalar multiple of a Haar unitary, and let $\nu_a$ be its Brown measure.  For $n\in\N$, there are constants $C(n)\asymp \sqrt{n}$ such that
\[ \|z^n\|_\infty = \sup_{\mathrm{supp}\,\nu_a}|z^n| \le C(n)\,\left[\int |z^n|^2\,d\nu_a(z,\bar{z})\right]^{1/2} = C(n)\,\|z^n\|_2. \]
\end{theorem}

\begin{proof} Fist note from Theorem \ref{theorem Brown measure of R-diag}, there is a function $f\colon[0,\|a\|_2]\to\R_+$ which is continuous and satisfies $f(\|a\|_2)>0$, such that $d\nu_a = f(r)dr\,d\theta$ with $\supp\,\nu_a$ equal to an annulus whose outer radius is $\|a\|_2$.  Of course, this means that $\sup_{\supp\,\nu_a} |z^n| = \|a\|_2^n$.  For the $2$-norm, let $M$ be the supremum of $f$ on $[0,\|a\|_2]$; then
\[ \int_{\supp\,\nu_a} |z^n|^2\,d\nu_a(z,\bar{z}) = \int_0^{2\pi}\int_0^{\|a\|_2} r^{2n}\,f(r)\,dr\,d\theta
\le \frac{2\pi M}{2n+1} \|a\|_2^{2n+1} = \frac{2\pi M\|a\|_2}{2n+1}\|z^n\|_\infty^2, \]
and this shows that $\|z^n\|_\infty/\|z^n\|_2 \apprge \sqrt{n}$.  For the reverse inequality, since $f$ is continuous and $f(\|a\|_2)>0$, there are $\e,m>0$ such that $f(r)\ge m>0$ for $r\in[\|a\|_2-\e,\|a\|_2]$, and so since $f\ge 0$ everywhere,
\[ \int_0^{2\pi}\int_0^{\|a\|_2} r^{2n}\,f(r)\,dr\,d\theta \ge 2\pi m \int_{\|a\|_2-\e}^{\|a\|_2} r^{2n}\,dr
= \frac{2\pi m\|a\|_2}{2n+1}\left(1-(1-\e/\|a\|_2)^{2n+1}\right)\|a\|_2^{2n} \apprge \frac{\|z^n\|_\infty^{2}}{n}. \]
\end{proof}

As discussed in Section \ref{subsection measures}, the Brown measure of a non-normal element $a$ (as most $\mathscr{R}$-diagonal elements are) does not respect mixed moments; that is, $\phi(a^\ast a) \ne \int |z|^2\,d\nu_a(z,\bar{z})$ in general, and so forth.  Nevertheless, as we see in Theorem \ref{theorem Brown measure Haagerup inequality}, a Haagerup inequality with the same $O(n^{1/2})$-behaviour holds in the space $\H L^2(\nu_a)$ of holomorphic $L^2$ functions with respect to the Brown measure of any $\mathscr{R}$-diagonal element.  $\H L^2(\nu_a)$ is, in some sense, the commutative model for our spaces $\H(a,I)$ (at least in the case where $|I|=1$), and so we see that the Brown measure does retain some information about mixed moments.

\section{Strong Ultracontractivity}\label{Strong Ultracontractivity}

In this final section, we apply our strong Haagerup inequality (Theorem \ref{main theorem}) to give strong ultracontractive bounds for the Ornstein-Uhlenbeck semigroup on $\H(a,I)$.  In Section \ref{subsection OU semigroups} we define said the O-U semigroup in this general context, and show that it is a natural generalization of the free O-U semigroup considered in \cite{Biane 2}.  In Section \ref{subsection ultracontractivity}, we prove optimal ultracontractive bounds, and discuss applications to free groups.

\subsection{Ornstein-Uhlenbeck semigroups}\label{subsection OU semigroups} Let $a$ be $\mathscr{R}$-diagonal.  Consider the operator $N_\mathrm{fin}$, defined on the algebraic direct sum $\bigoplus_{n=0}^\infty \H^{(n)}(a,I)$ (which is, of course, dense in $L^2(\H(a,I),\phi)$) as the linear extension of $N_\mathrm{fin}(h_n) = n\,h_n$ for $h_n\in\H^{(n)}(a,I)$.  Since $h_n\perp h_m$ for $n\ne m$ (this follows from the $\ast$-freeness of the $a_i$), the operator $N_\mathrm{fin}$ is symmetric and lower-semi-bounded by $0$.  Thus, by the Friedrich's extension theorem, $N_\mathrm{fin}$ extends to a densely-defined (unbounded) self-adjoint operator $N$ on $L^2(\H(a,I),\phi)$, and this operator is postive semidefinite.  We will refer to $N$ as the {\bf number operator} affiliated with $\H(a,I)$.

\begin{proposition}\label{prop semigroup} The number operator $N$ affiliated with $\H(a,I)$ generates a $\mathscr{C}_0$ contraction semigroup $e^{-tN}$ on $L^2(\H(a,I),\phi)$.
\end{proposition}

\begin{proof} Since the spaces $\H^{(n)}(a,I)$ reduce $N$, we see easily that $e^{-tN}$ must act via
\[ e^{-tN}\sum_{n=0}^\infty h_n = \sum_{n=0}^\infty e^{-nt} h_n. \]
It is then immediately verified that $e^{-tN}$ is a contraction semigroup, since $e^{-nt}\le 1$ for all $t\ge 0$.  To prove that is it $\mathscr{C}_0$, it suffices to show that $w$-$\lim_{t\downarrow 0}e^{-tN}h = h$ for each $h\in L^2(\H(a,I),\phi)$.  Let $h=\sum h_n$ and $g=\sum g_n$; since $h_n\perp g_m$ for $n\ne m$,
\[ \langle e^{-tN}h,g\rangle = \left\langle \sum_{n=0}^\infty e^{-nt}h_n,\sum_{m=0}^\infty g_m\right\rangle
=\sum_{n=0}^\infty e^{-nt}\langle h_n,g_n\rangle. \]
As both $h$ and $g$ are in $L^2$, the sequence $\langle h_n,g_n\rangle$ is in $\ell^1$, and since $e^{-nt}\le 1$, it follows from the dominated convergence theorem that
\[ \lim_{t\downarrow 0} \sum_{n=0}^\infty e^{-nt}\langle h_n,g_n\rangle = \sum_{n=0}^\infty \langle h_n,g_n \rangle = \langle h,g \rangle. \]  \end{proof}

\medskip

An important example of this number operator is given in the case of a circular element $a=c$.  In this case, $\H(c,I)$ is naturally isomorphic to the holomorphic space $\H_0(\HH)$ over a Hilbert space $\HH$ of dimension $|I|$, as defined in the first author's paper \cite{Kemp}, and the number operator $N$ above is just the free Ornstein-Uhlenbeck (number) operator $N_0$ considered in that paper.  $N_0$ is the restriction to the holomorphic space $\H_0(\HH)$ of the free Ornstein-Uhlenbeck operator defined in \cite{Biane 2} on the free group factor $L(\F_{|I|})$, which coincides with the $0$-Gaussian factor $\Gamma_0(\HH)$ introduced in \cite{Voiculescu} and further developed in \cite{Bozejko Speicher, BoKS}.  There is a family of such spaces $\Gamma_q(\HH)$ for $-1\le q\le 1$ (with $q=1$ corresponding to the classical theory of Gaussian random variables, and $q=-1$ the hyperfinite $\mathrm{II}_1$-factor), and Biane introduced number operators $N_q$ affiliated to each of them.  We should also note that, in \cite{Biane 1}, Biane introduced a space isomorphic to $\H(c,I)$, but did not consider the action of a number operator on it.

\medskip

The main theorem of \cite{Kemp} shows as a special case (the case $q=0$) that the semigroup $e^{-tN}$ affiliated with $\H(c,I)$ is not only a contraction semigroup on $L^2(\H(c,I),\phi)$ (for tracial $\phi$), but is in fact {\em strongly hypercontractive}:

\begin{theorem}[Theorem 4 in \cite{Kemp}] Let $r>2$ be an even integer, and let $t_J(2,r)=\frac{1}{2}\log\frac{r}{2}$.  Then for $t\ge t_J(2,r)$, $e^{-tN}$ is a contraction from $L^2(\H(c,I),\phi)$ to $L^r(\H(c,I),\phi)$. \end{theorem}

This {\em strong hypercontractivity} theorem is the precise analogue of the same theorem in the context of the spaces $\H L^r(\C^n,\gamma)$ (where $\gamma$ is Gauss measure) proved by Janson in \cite{Janson}.  (We should note, however, that Janson's theorem holds from $L^p\to L^r$ for $0<p\le r<\infty$, not just the discrete values in \cite{Kemp}.)  The time $t_J$ is shorter than the least time to contraction $t_N$ in the real spaces $L^r(\R^n,\gamma)$, where the hypercontractivity inequalities were first proved and studied by Nelson in \cite{Nelson}.  The main theorem of \cite{Biane 2} is the generalization of Nelson's hypercontractivity theorem to the $q$-Gaussian factors.

\subsection{Ultracontractivity}\label{subsection ultracontractivity} In the classical holomorphic case studied by Janson, while the semigroup $e^{-tN}$ is a contractive map from $\H L^2$ to $\H L^r$ for any $r>2$, once $t$ is large enough it is also unbounded for $t < t_J(2,r)$.  As a result, the semigroup $e^{-tN}$ does not map $\H L^2$ into the algebra of bounded functions for any time.  Of course, in the classical context, the algebra of bounded functions contains no holomorphic functions save constants; even in the full real spaces, the same effect holds.  This is essentially due to the fact that the kernel of the semigroup $e^{-tN}$ in these cases, the {\em Mehler kernel}, is not a bounded function.

\medskip

A semigroup is called {\em ultracontractive} if it maps $L^2$ into $L^\infty$ for all $t>0$.  The Ornstein-Uhlenbeck semigroups studied by Nelson and Janson (and many others) fail to be ultracontractive.  Nevertheless, the non-commutative counterpart $e^{-tN_0}$ on the free group factor {\em is} ultracontractive, as shown in \cite{Biane 2} and essentially in \cite{Bozejko 1}.

\begin{proposition}[Corollary 3 in \cite{Biane 2}]\label{prop Biane ultracontractivity} The free Ornstein-Uhlenbeck semigroup $e^{-tN_0}$ is ultracontractive; there is $c>0$ with
\[ \|e^{-tN_0}\|_{L^2(\Gamma_0)\to L^\infty(\Gamma_0)} \le c\, t^{-3/2}, \quad 0<t<1. \]
\end{proposition}
(In general the function $t\mapsto \|e^{-tN_0}X\|_r$ is decreasing for any $X$ and $r$, hence it is only small-time behaviour which is interesting.)  Bo\.zejko later generalized this theorem to all the $\Gamma_q$ factors with $-1<q<1$; see \cite{Bozejko 2}.

\medskip

The generators of the algebra $\Gamma_0$ (the free group factor) are $\ast$-free semicircular elements.  Thus, the $\ast$-algebra generated by $\H(c,I)$ is contained in $\Gamma_0$, and the ultracontractive $O(t^{-3/2})$-bound of Proposition \ref{prop Biane ultracontractivity} also holds for the semigroup $e^{-tN}$ affiliated with $\H(c,I)$ defined above.  Using our main theorem, Theorem \ref{main theorem}, we may essentially follow Biane's argument and prove a stronger form of Proposition \ref{prop Biane ultracontractivity} not only for the algebra $\H(c,I)\cong \H_0(\HH)$, but in fact for all $\H(a,I)$ with $a$ $\mathscr{R}$-diagonal.  Indeed, we find that the short-time behaviour in the $\mathscr{R}$-diagonal case is $O(t^{-1})$.

\begin{theorem}\label{theorem strong ultracontractivity} Let $a$ be $\mathscr{R}$-diagonal, and let $N$ be the number operator affiliated with $\H(a,I)$.  Then $e^{-tN}$ is ultracontractive; for each $h\in L^2(\H(a,I),\phi)$, $e^{-tN}h\in\H(a,I)$ for $t>0$, and moreover
\begin{equation}\label{strong ultracontractivity bound} \|e^{-tN}h\| \le \frac{1}{2}C_a\,t^{-1}\,\|h\|_2 \quad t>0. \end{equation}
\end{theorem}
(Here $C_a$ is the same constant as in Theorem \ref{main theorem}.) We refer to Theorem \ref{theorem strong ultracontractivity} as {\em strong ultracontractivity}, as it is a stronger version of the inequality in Proposition \ref{prop Biane ultracontractivity} which holds when the semigroup is restricted to a holomorphic subspace.  This is similar in spirit to the stronger form of hypercontractivity \cite{Janson} which holds in the holomorphic version of Nelson's setup in \cite{Nelson}.  We emphasize, again, that ultracontractivity is a {\em strictly} non-commutative effect in this case, since the semigroup is unbounded from $L^2\to L^\infty$ in the classical (real and holomorphic) contexts.  Theorem \ref{theorem strong ultracontractivity} is thus an essentially non-commutative result which highlights the interesting phenomenon that many functional inequalities improve in the holomorphic category.

\begin{proof} Let $h = \sum_{n=0}^\infty h_n$ with $h_n\in\H^{(n)}(a,I)$. We estimate
\[ \|e^{-tN}h\| = \left\|\sum_{n=0}^\infty e^{-nt}h_n\right\| \le \sum_{n=0}^\infty e^{-nt}\|h_n\|. \]
We now employ Theorem \ref{main theorem}, which implies that $h_n\in\H(a,I)$ and $\|h_n\|\le C_a\,\sqrt{n}\,\|h_n\|_2$.  Thus,
\[ \|e^{-tN}h\| \le C_a\sum_{n=0}^\infty \sqrt{n}\,e^{-nt}\|h_n\|_2 \le C_a\left[\sum_{n=0}^\infty n\,e^{-2nt}\right]^{1/2}\cdot\left[\sum_{n=0}^\infty \|h_n\|_2^2\right]^{1/2}, \]
where we have used the Cauchy-Schwarz inequality.  The second factor is just $\|h\|_2^2$.  The first factor is the derivative of $-\frac{1}{2}\sum_{n=0}^\infty e^{-2nt} = -\frac{1}{2}\frac{1}{1-e^{-2t}}$.  The reader may readily verify that we thus have
\[ \|e^{-tN}h\| \le C_a \frac{e^{-t}}{1-e^{-2t}}\,\|h\|_2 \]
for all $t>0$.  This shows that $e^{-tN}h\in\H(a,I)$.  Moreover, the function $t\mapsto \frac{te^{-t}}{1-e^{-2t}}$ is decreasing on $\R_+$ and has limit $1/2$ at $t=0$.  This proves Equation \ref{strong ultracontractivity bound}.
\end{proof}

It is typical to prove, from a bound like Equation \ref{strong ultracontractivity bound}, a Sobolev inequality of the form $\|h\|_p \le c\langle Nh,h\rangle, h\in\mathscr{D}(N)$ for an appropriate $p>2$; indeed, if $e^{-tN}$ in Theorem \ref{theorem strong ultracontractivity} were a classical sub-Markovian semigroup defined on $L^2$ of a Radon measure, we could use the standard techniques in, for example, \cite{CSV}, to prove a strong Sobolev imbedding theorem (for {\em any} $p<\infty$) in this case.  However, the techniques necessary to implement such a proof use the Marcinkewicz interpolation theorem in a fundamental way.  As pointed out in \cite{Kemp}, holomorphic spaces like $\H(a,I)$ (in particular in the case $a=c$) tend not to be complex interpolation scale (at least in the $|I|=\infty$ case).  Thus, we are unable to prove a Sobolev inequality for $\H(a,I)$ using known-techniques.

\medskip

We finally remark that one interesting new application of this theorem is to the discrete O-U semigroup on the free semigroup $\F_k$ (or rather its restriction to $\F_k^+$).  As noted above, the algebra $\H(u,I_k)$ with $u$ a Haar unitary and $|I_k|=k$ is isomorphic to the convolution-norm closure of $\F_k^+$ in $L(\F_k)$, and thus $L^2(\H(u,I_k),\phi_k)\cong \ell^2(\F_k^+)$, where the number operator $N$ acts by $Nw = nw$ on a word $w$ of length $n$.  The same semigroup $e^{-tN}$ defined on all of $\F_k$ was essentially introduced in \cite{Haagerup}, and has been studied in \cite{JLX, JX} with a view towards $L^p$-contraction bounds; to the authors' knowledge, Theorem \ref{theorem strong ultracontractivity} yields the first ultracontractive bound in that context.

\medskip

\noindent {\bf Acknowledgement.}  The authors kindly thank Drew Armstrong and Heather Armstrong for their combinatorial insights.

\end{document}

%% file: fig1.pstex_t
\begin{picture}(0,0)%
\includegraphics{fig1.pstex}%
\end{picture}%
\setlength{\unitlength}{3947sp}%
\begingroup\makeatletter\ifx\SetFigFont\undefined%
\gdef\SetFigFont#1#2#3#4#5{%
  \reset@font\fontsize{#1}{#2pt}%
  \fontfamily{#3}\fontseries{#4}\fontshape{#5}%
  \selectfont}%
\fi\endgroup%
\begin{picture}(5307,2598)(1156,-2819)
\put(1156,-353){\makebox(0,0)[lb]{\smash{\SetFigFont{12}{14.4}{\familydefault}{\mddefault}{\updefault}{\color[rgb]{0,0,0}$1$}%
}}}
\put(1606,-353){\makebox(0,0)[lb]{\smash{\SetFigFont{12}{14.4}{\familydefault}{\mddefault}{\updefault}{\color[rgb]{0,0,0}$2$}%
}}}
\put(2056,-353){\makebox(0,0)[lb]{\smash{\SetFigFont{12}{14.4}{\familydefault}{\mddefault}{\updefault}{\color[rgb]{0,0,0}$3$}%
}}}
\put(2506,-353){\makebox(0,0)[lb]{\smash{\SetFigFont{12}{14.4}{\familydefault}{\mddefault}{\updefault}{\color[rgb]{0,0,0}$4$}%
}}}
\put(2956,-353){\makebox(0,0)[lb]{\smash{\SetFigFont{12}{14.4}{\familydefault}{\mddefault}{\updefault}{\color[rgb]{0,0,0}$5$}%
}}}
\put(3406,-353){\makebox(0,0)[lb]{\smash{\SetFigFont{12}{14.4}{\familydefault}{\mddefault}{\updefault}{\color[rgb]{0,0,0}$6$}%
}}}
\put(4156,-353){\makebox(0,0)[lb]{\smash{\SetFigFont{12}{14.4}{\familydefault}{\mddefault}{\updefault}{\color[rgb]{0,0,0}$1$}%
}}}
\put(4606,-353){\makebox(0,0)[lb]{\smash{\SetFigFont{12}{14.4}{\familydefault}{\mddefault}{\updefault}{\color[rgb]{0,0,0}$2$}%
}}}
\put(5056,-353){\makebox(0,0)[lb]{\smash{\SetFigFont{12}{14.4}{\familydefault}{\mddefault}{\updefault}{\color[rgb]{0,0,0}$3$}%
}}}
\put(5506,-353){\makebox(0,0)[lb]{\smash{\SetFigFont{12}{14.4}{\familydefault}{\mddefault}{\updefault}{\color[rgb]{0,0,0}$4$}%
}}}
\put(5956,-353){\makebox(0,0)[lb]{\smash{\SetFigFont{12}{14.4}{\familydefault}{\mddefault}{\updefault}{\color[rgb]{0,0,0}$5$}%
}}}
\put(6406,-353){\makebox(0,0)[lb]{\smash{\SetFigFont{12}{14.4}{\familydefault}{\mddefault}{\updefault}{\color[rgb]{0,0,0}$6$}%
}}}
\put(1156,-1853){\makebox(0,0)[lb]{\smash{\SetFigFont{12}{14.4}{\familydefault}{\mddefault}{\updefault}{\color[rgb]{0,0,0}$1$}%
}}}
\put(1606,-1853){\makebox(0,0)[lb]{\smash{\SetFigFont{12}{14.4}{\familydefault}{\mddefault}{\updefault}{\color[rgb]{0,0,0}$2$}%
}}}
\put(2056,-1853){\makebox(0,0)[lb]{\smash{\SetFigFont{12}{14.4}{\familydefault}{\mddefault}{\updefault}{\color[rgb]{0,0,0}$3$}%
}}}
\put(2506,-1853){\makebox(0,0)[lb]{\smash{\SetFigFont{12}{14.4}{\familydefault}{\mddefault}{\updefault}{\color[rgb]{0,0,0}$4$}%
}}}
\put(2956,-1853){\makebox(0,0)[lb]{\smash{\SetFigFont{12}{14.4}{\familydefault}{\mddefault}{\updefault}{\color[rgb]{0,0,0}$5$}%
}}}
\put(3406,-1853){\makebox(0,0)[lb]{\smash{\SetFigFont{12}{14.4}{\familydefault}{\mddefault}{\updefault}{\color[rgb]{0,0,0}$6$}%
}}}
\put(4156,-1853){\makebox(0,0)[lb]{\smash{\SetFigFont{12}{14.4}{\familydefault}{\mddefault}{\updefault}{\color[rgb]{0,0,0}$1$}%
}}}
\put(4606,-1853){\makebox(0,0)[lb]{\smash{\SetFigFont{12}{14.4}{\familydefault}{\mddefault}{\updefault}{\color[rgb]{0,0,0}$2$}%
}}}
\put(5056,-1853){\makebox(0,0)[lb]{\smash{\SetFigFont{12}{14.4}{\familydefault}{\mddefault}{\updefault}{\color[rgb]{0,0,0}$3$}%
}}}
\put(5506,-1853){\makebox(0,0)[lb]{\smash{\SetFigFont{12}{14.4}{\familydefault}{\mddefault}{\updefault}{\color[rgb]{0,0,0}$4$}%
}}}
\put(5956,-1853){\makebox(0,0)[lb]{\smash{\SetFigFont{12}{14.4}{\familydefault}{\mddefault}{\updefault}{\color[rgb]{0,0,0}$5$}%
}}}
\put(6406,-1853){\makebox(0,0)[lb]{\smash{\SetFigFont{12}{14.4}{\familydefault}{\mddefault}{\updefault}{\color[rgb]{0,0,0}$6$}%
}}}
\put(2251,-1261){\makebox(0,0)[lb]{\smash{\SetFigFont{12}{14.4}{\familydefault}{\mddefault}{\updefault}{\color[rgb]{0,0,0}$1_6$}%
}}}
\put(5251,-2761){\makebox(0,0)[lb]{\smash{\SetFigFont{12}{14.4}{\familydefault}{\mddefault}{\updefault}{\color[rgb]{0,0,0}$0_6$}%
}}}
\end{picture}

%% file: fig2.pstex_t
\begin{picture}(0,0)%
\includegraphics{fig2.pstex}%
\end{picture}%
\setlength{\unitlength}{3947sp}%
\begingroup\makeatletter\ifx\SetFigFont\undefined%
\gdef\SetFigFont#1#2#3#4#5{%
  \reset@font\fontsize{#1}{#2pt}%
  \fontfamily{#3}\fontseries{#4}\fontshape{#5}%
  \selectfont}%
\fi\endgroup%
\begin{picture}(6960,3529)(1201,-3158)
\put(1201,239){\makebox(0,0)[lb]{\smash{\SetFigFont{12}{14.4}{\familydefault}{\mddefault}{\updefault}{\color[rgb]{0,0,0}$c$}%
}}}
\put(1501,239){\makebox(0,0)[lb]{\smash{\SetFigFont{12}{14.4}{\familydefault}{\mddefault}{\updefault}{\color[rgb]{0,0,0}$c$}%
}}}
\put(1801,239){\makebox(0,0)[lb]{\smash{\SetFigFont{12}{14.4}{\familydefault}{\mddefault}{\updefault}{\color[rgb]{0,0,0}$c$}%
}}}
\put(2101,239){\makebox(0,0)[lb]{\smash{\SetFigFont{12}{14.4}{\familydefault}{\mddefault}{\updefault}{\color[rgb]{0,0,0}$c^\ast$}%
}}}
\put(2401,239){\makebox(0,0)[lb]{\smash{\SetFigFont{12}{14.4}{\familydefault}{\mddefault}{\updefault}{\color[rgb]{0,0,0}$c^\ast$}%
}}}
\put(2701,239){\makebox(0,0)[lb]{\smash{\SetFigFont{12}{14.4}{\familydefault}{\mddefault}{\updefault}{\color[rgb]{0,0,0}$c^\ast$}%
}}}
\put(3001,239){\makebox(0,0)[lb]{\smash{\SetFigFont{12}{14.4}{\familydefault}{\mddefault}{\updefault}{\color[rgb]{0,0,0}$c$}%
}}}
\put(3301,239){\makebox(0,0)[lb]{\smash{\SetFigFont{12}{14.4}{\familydefault}{\mddefault}{\updefault}{\color[rgb]{0,0,0}$c$}%
}}}
\put(3601,239){\makebox(0,0)[lb]{\smash{\SetFigFont{12}{14.4}{\familydefault}{\mddefault}{\updefault}{\color[rgb]{0,0,0}$c$}%
}}}
\put(3901,239){\makebox(0,0)[lb]{\smash{\SetFigFont{12}{14.4}{\familydefault}{\mddefault}{\updefault}{\color[rgb]{0,0,0}$c^\ast$}%
}}}
\put(4201,239){\makebox(0,0)[lb]{\smash{\SetFigFont{12}{14.4}{\familydefault}{\mddefault}{\updefault}{\color[rgb]{0,0,0}$c^\ast$}%
}}}
\put(4501,239){\makebox(0,0)[lb]{\smash{\SetFigFont{12}{14.4}{\familydefault}{\mddefault}{\updefault}{\color[rgb]{0,0,0}$c^\ast$}%
}}}
\put(4801,239){\makebox(0,0)[lb]{\smash{\SetFigFont{12}{14.4}{\familydefault}{\mddefault}{\updefault}{\color[rgb]{0,0,0}$c$}%
}}}
\put(5101,239){\makebox(0,0)[lb]{\smash{\SetFigFont{12}{14.4}{\familydefault}{\mddefault}{\updefault}{\color[rgb]{0,0,0}$c$}%
}}}
\put(5401,239){\makebox(0,0)[lb]{\smash{\SetFigFont{12}{14.4}{\familydefault}{\mddefault}{\updefault}{\color[rgb]{0,0,0}$c$}%
}}}
\put(5701,239){\makebox(0,0)[lb]{\smash{\SetFigFont{12}{14.4}{\familydefault}{\mddefault}{\updefault}{\color[rgb]{0,0,0}$c^\ast$}%
}}}
\put(6001,239){\makebox(0,0)[lb]{\smash{\SetFigFont{12}{14.4}{\familydefault}{\mddefault}{\updefault}{\color[rgb]{0,0,0}$c^\ast$}%
}}}
\put(6301,239){\makebox(0,0)[lb]{\smash{\SetFigFont{12}{14.4}{\familydefault}{\mddefault}{\updefault}{\color[rgb]{0,0,0}$c^\ast$}%
}}}
\put(6601,239){\makebox(0,0)[lb]{\smash{\SetFigFont{12}{14.4}{\familydefault}{\mddefault}{\updefault}{\color[rgb]{0,0,0}$c$}%
}}}
\put(6901,239){\makebox(0,0)[lb]{\smash{\SetFigFont{12}{14.4}{\familydefault}{\mddefault}{\updefault}{\color[rgb]{0,0,0}$c$}%
}}}
\put(7201,239){\makebox(0,0)[lb]{\smash{\SetFigFont{12}{14.4}{\familydefault}{\mddefault}{\updefault}{\color[rgb]{0,0,0}$c$}%
}}}
\put(7501,239){\makebox(0,0)[lb]{\smash{\SetFigFont{12}{14.4}{\familydefault}{\mddefault}{\updefault}{\color[rgb]{0,0,0}$c^\ast$}%
}}}
\put(7801,239){\makebox(0,0)[lb]{\smash{\SetFigFont{12}{14.4}{\familydefault}{\mddefault}{\updefault}{\color[rgb]{0,0,0}$c^\ast$}%
}}}
\put(8101,239){\makebox(0,0)[lb]{\smash{\SetFigFont{12}{14.4}{\familydefault}{\mddefault}{\updefault}{\color[rgb]{0,0,0}$c^\ast$}%
}}}
\put(1201,-1711){\makebox(0,0)[lb]{\smash{\SetFigFont{12}{14.4}{\familydefault}{\mddefault}{\updefault}{\color[rgb]{0,0,0}$c$}%
}}}
\put(1501,-1711){\makebox(0,0)[lb]{\smash{\SetFigFont{12}{14.4}{\familydefault}{\mddefault}{\updefault}{\color[rgb]{0,0,0}$c$}%
}}}
\put(1801,-1711){\makebox(0,0)[lb]{\smash{\SetFigFont{12}{14.4}{\familydefault}{\mddefault}{\updefault}{\color[rgb]{0,0,0}$c$}%
}}}
\put(2101,-1711){\makebox(0,0)[lb]{\smash{\SetFigFont{12}{14.4}{\familydefault}{\mddefault}{\updefault}{\color[rgb]{0,0,0}$c^\ast$}%
}}}
\put(2401,-1711){\makebox(0,0)[lb]{\smash{\SetFigFont{12}{14.4}{\familydefault}{\mddefault}{\updefault}{\color[rgb]{0,0,0}$c^\ast$}%
}}}
\put(2701,-1711){\makebox(0,0)[lb]{\smash{\SetFigFont{12}{14.4}{\familydefault}{\mddefault}{\updefault}{\color[rgb]{0,0,0}$c^\ast$}%
}}}
\put(3001,-1711){\makebox(0,0)[lb]{\smash{\SetFigFont{12}{14.4}{\familydefault}{\mddefault}{\updefault}{\color[rgb]{0,0,0}$c$}%
}}}
\put(3301,-1711){\makebox(0,0)[lb]{\smash{\SetFigFont{12}{14.4}{\familydefault}{\mddefault}{\updefault}{\color[rgb]{0,0,0}$c$}%
}}}
\put(3601,-1711){\makebox(0,0)[lb]{\smash{\SetFigFont{12}{14.4}{\familydefault}{\mddefault}{\updefault}{\color[rgb]{0,0,0}$c$}%
}}}
\put(3901,-1711){\makebox(0,0)[lb]{\smash{\SetFigFont{12}{14.4}{\familydefault}{\mddefault}{\updefault}{\color[rgb]{0,0,0}$c^\ast$}%
}}}
\put(4201,-1711){\makebox(0,0)[lb]{\smash{\SetFigFont{12}{14.4}{\familydefault}{\mddefault}{\updefault}{\color[rgb]{0,0,0}$c^\ast$}%
}}}
\put(4501,-1711){\makebox(0,0)[lb]{\smash{\SetFigFont{12}{14.4}{\familydefault}{\mddefault}{\updefault}{\color[rgb]{0,0,0}$c^\ast$}%
}}}
\put(4801,-1711){\makebox(0,0)[lb]{\smash{\SetFigFont{12}{14.4}{\familydefault}{\mddefault}{\updefault}{\color[rgb]{0,0,0}$c$}%
}}}
\put(5101,-1711){\makebox(0,0)[lb]{\smash{\SetFigFont{12}{14.4}{\familydefault}{\mddefault}{\updefault}{\color[rgb]{0,0,0}$c$}%
}}}
\put(5401,-1711){\makebox(0,0)[lb]{\smash{\SetFigFont{12}{14.4}{\familydefault}{\mddefault}{\updefault}{\color[rgb]{0,0,0}$c$}%
}}}
\put(5701,-1711){\makebox(0,0)[lb]{\smash{\SetFigFont{12}{14.4}{\familydefault}{\mddefault}{\updefault}{\color[rgb]{0,0,0}$c^\ast$}%
}}}
\put(6001,-1711){\makebox(0,0)[lb]{\smash{\SetFigFont{12}{14.4}{\familydefault}{\mddefault}{\updefault}{\color[rgb]{0,0,0}$c^\ast$}%
}}}
\put(6301,-1711){\makebox(0,0)[lb]{\smash{\SetFigFont{12}{14.4}{\familydefault}{\mddefault}{\updefault}{\color[rgb]{0,0,0}$c^\ast$}%
}}}
\put(6601,-1711){\makebox(0,0)[lb]{\smash{\SetFigFont{12}{14.4}{\familydefault}{\mddefault}{\updefault}{\color[rgb]{0,0,0}$c$}%
}}}
\put(6901,-1711){\makebox(0,0)[lb]{\smash{\SetFigFont{12}{14.4}{\familydefault}{\mddefault}{\updefault}{\color[rgb]{0,0,0}$c$}%
}}}
\put(7201,-1711){\makebox(0,0)[lb]{\smash{\SetFigFont{12}{14.4}{\familydefault}{\mddefault}{\updefault}{\color[rgb]{0,0,0}$c$}%
}}}
\put(7501,-1711){\makebox(0,0)[lb]{\smash{\SetFigFont{12}{14.4}{\familydefault}{\mddefault}{\updefault}{\color[rgb]{0,0,0}$c^\ast$}%
}}}
\put(7801,-1711){\makebox(0,0)[lb]{\smash{\SetFigFont{12}{14.4}{\familydefault}{\mddefault}{\updefault}{\color[rgb]{0,0,0}$c^\ast$}%
}}}
\put(8101,-1711){\makebox(0,0)[lb]{\smash{\SetFigFont{12}{14.4}{\familydefault}{\mddefault}{\updefault}{\color[rgb]{0,0,0}$c^\ast$}%
}}}
\end{picture}

%% file: fig3.pstex_t
\begin{picture}(0,0)%
\includegraphics{fig3.pstex}%
\end{picture}%
\setlength{\unitlength}{2960sp}%
\begingroup\makeatletter\ifx\SetFigFont\undefined%
\gdef\SetFigFont#1#2#3#4#5{%
  \reset@font\fontsize{#1}{#2pt}%
  \fontfamily{#3}\fontseries{#4}\fontshape{#5}%
  \selectfont}%
\fi\endgroup%
\begin{picture}(8467,4007)(1626,-4490)
\put(3351,-616){\makebox(0,0)[lb]{\smash{\SetFigFont{9}{10.8}{\familydefault}{\mddefault}{\updefault}{\color[rgb]{0,0,0}$c^\ast$}%
}}}
\put(3861,-631){\makebox(0,0)[lb]{\smash{\SetFigFont{9}{10.8}{\familydefault}{\mddefault}{\updefault}{\color[rgb]{0,0,0}$c$}%
}}}
\put(4336,-766){\makebox(0,0)[lb]{\smash{\SetFigFont{9}{10.8}{\familydefault}{\mddefault}{\updefault}{\color[rgb]{0,0,0}$c$}%
}}}
\put(4761,-1026){\makebox(0,0)[lb]{\smash{\SetFigFont{9}{10.8}{\familydefault}{\mddefault}{\updefault}{\color[rgb]{0,0,0}$c$}%
}}}
\put(5106,-1381){\makebox(0,0)[lb]{\smash{\SetFigFont{9}{10.8}{\familydefault}{\mddefault}{\updefault}{\color[rgb]{0,0,0}$c^\ast$}%
}}}
\put(5356,-1826){\makebox(0,0)[lb]{\smash{\SetFigFont{9}{10.8}{\familydefault}{\mddefault}{\updefault}{\color[rgb]{0,0,0}$c^\ast$}%
}}}
\put(5476,-2306){\makebox(0,0)[lb]{\smash{\SetFigFont{9}{10.8}{\familydefault}{\mddefault}{\updefault}{\color[rgb]{0,0,0}$c^\ast$}%
}}}
\put(5466,-2821){\makebox(0,0)[lb]{\smash{\SetFigFont{9}{10.8}{\familydefault}{\mddefault}{\updefault}{\color[rgb]{0,0,0}$c$}%
}}}
\put(5326,-3296){\makebox(0,0)[lb]{\smash{\SetFigFont{9}{10.8}{\familydefault}{\mddefault}{\updefault}{\color[rgb]{0,0,0}$c$}%
}}}
\put(5056,-3731){\makebox(0,0)[lb]{\smash{\SetFigFont{9}{10.8}{\familydefault}{\mddefault}{\updefault}{\color[rgb]{0,0,0}$c$}%
}}}
\put(4706,-4066){\makebox(0,0)[lb]{\smash{\SetFigFont{9}{10.8}{\familydefault}{\mddefault}{\updefault}{\color[rgb]{0,0,0}$c^\ast$}%
}}}
\put(4276,-4311){\makebox(0,0)[lb]{\smash{\SetFigFont{9}{10.8}{\familydefault}{\mddefault}{\updefault}{\color[rgb]{0,0,0}$c^\ast$}%
}}}
\put(3776,-4441){\makebox(0,0)[lb]{\smash{\SetFigFont{9}{10.8}{\familydefault}{\mddefault}{\updefault}{\color[rgb]{0,0,0}$c^\ast$}%
}}}
\put(3266,-4426){\makebox(0,0)[lb]{\smash{\SetFigFont{9}{10.8}{\familydefault}{\mddefault}{\updefault}{\color[rgb]{0,0,0}$c$}%
}}}
\put(2781,-4286){\makebox(0,0)[lb]{\smash{\SetFigFont{9}{10.8}{\familydefault}{\mddefault}{\updefault}{\color[rgb]{0,0,0}$c$}%
}}}
\put(1656,-2231){\makebox(0,0)[lb]{\smash{\SetFigFont{9}{10.8}{\familydefault}{\mddefault}{\updefault}{\color[rgb]{0,0,0}$c$}%
}}}
\put(1801,-1746){\makebox(0,0)[lb]{\smash{\SetFigFont{9}{10.8}{\familydefault}{\mddefault}{\updefault}{\color[rgb]{0,0,0}$c$}%
}}}
\put(2061,-1326){\makebox(0,0)[lb]{\smash{\SetFigFont{9}{10.8}{\familydefault}{\mddefault}{\updefault}{\color[rgb]{0,0,0}$c$}%
}}}
\put(2416,-986){\makebox(0,0)[lb]{\smash{\SetFigFont{9}{10.8}{\familydefault}{\mddefault}{\updefault}{\color[rgb]{0,0,0}$c^\ast$}%
}}}
\put(2861,-741){\makebox(0,0)[lb]{\smash{\SetFigFont{9}{10.8}{\familydefault}{\mddefault}{\updefault}{\color[rgb]{0,0,0}$c^\ast$}%
}}}
\put(2356,-4026){\makebox(0,0)[lb]{\smash{\SetFigFont{9}{10.8}{\familydefault}{\mddefault}{\updefault}{\color[rgb]{0,0,0}$c$}%
}}}
\put(1961,-3683){\makebox(0,0)[lb]{\smash{\SetFigFont{9}{10.8}{\familydefault}{\mddefault}{\updefault}{\color[rgb]{0,0,0}$c^\ast$}%
}}}
\put(1719,-3235){\makebox(0,0)[lb]{\smash{\SetFigFont{9}{10.8}{\familydefault}{\mddefault}{\updefault}{\color[rgb]{0,0,0}$c^\ast$}%
}}}
\put(1626,-2744){\makebox(0,0)[lb]{\smash{\SetFigFont{9}{10.8}{\familydefault}{\mddefault}{\updefault}{\color[rgb]{0,0,0}$c^\ast$}%
}}}
\put(7957,-615){\makebox(0,0)[lb]{\smash{\SetFigFont{9}{10.8}{\familydefault}{\mddefault}{\updefault}{\color[rgb]{0,0,0}$c^\ast$}%
}}}
\put(8467,-630){\makebox(0,0)[lb]{\smash{\SetFigFont{9}{10.8}{\familydefault}{\mddefault}{\updefault}{\color[rgb]{0,0,0}$c$}%
}}}
\put(8942,-765){\makebox(0,0)[lb]{\smash{\SetFigFont{9}{10.8}{\familydefault}{\mddefault}{\updefault}{\color[rgb]{0,0,0}$c$}%
}}}
\put(9367,-1025){\makebox(0,0)[lb]{\smash{\SetFigFont{9}{10.8}{\familydefault}{\mddefault}{\updefault}{\color[rgb]{0,0,0}$c$}%
}}}
\put(9712,-1380){\makebox(0,0)[lb]{\smash{\SetFigFont{9}{10.8}{\familydefault}{\mddefault}{\updefault}{\color[rgb]{0,0,0}$c^\ast$}%
}}}
\put(9962,-1825){\makebox(0,0)[lb]{\smash{\SetFigFont{9}{10.8}{\familydefault}{\mddefault}{\updefault}{\color[rgb]{0,0,0}$c^\ast$}%
}}}
\put(10082,-2305){\makebox(0,0)[lb]{\smash{\SetFigFont{9}{10.8}{\familydefault}{\mddefault}{\updefault}{\color[rgb]{0,0,0}$c^\ast$}%
}}}
\put(10072,-2820){\makebox(0,0)[lb]{\smash{\SetFigFont{9}{10.8}{\familydefault}{\mddefault}{\updefault}{\color[rgb]{0,0,0}$c$}%
}}}
\put(9932,-3295){\makebox(0,0)[lb]{\smash{\SetFigFont{9}{10.8}{\familydefault}{\mddefault}{\updefault}{\color[rgb]{0,0,0}$c$}%
}}}
\put(9662,-3730){\makebox(0,0)[lb]{\smash{\SetFigFont{9}{10.8}{\familydefault}{\mddefault}{\updefault}{\color[rgb]{0,0,0}$c$}%
}}}
\put(9312,-4065){\makebox(0,0)[lb]{\smash{\SetFigFont{9}{10.8}{\familydefault}{\mddefault}{\updefault}{\color[rgb]{0,0,0}$c^\ast$}%
}}}
\put(8882,-4310){\makebox(0,0)[lb]{\smash{\SetFigFont{9}{10.8}{\familydefault}{\mddefault}{\updefault}{\color[rgb]{0,0,0}$c^\ast$}%
}}}
\put(8382,-4440){\makebox(0,0)[lb]{\smash{\SetFigFont{9}{10.8}{\familydefault}{\mddefault}{\updefault}{\color[rgb]{0,0,0}$c^\ast$}%
}}}
\put(7872,-4425){\makebox(0,0)[lb]{\smash{\SetFigFont{9}{10.8}{\familydefault}{\mddefault}{\updefault}{\color[rgb]{0,0,0}$c$}%
}}}
\put(7387,-4285){\makebox(0,0)[lb]{\smash{\SetFigFont{9}{10.8}{\familydefault}{\mddefault}{\updefault}{\color[rgb]{0,0,0}$c$}%
}}}
\put(6262,-2230){\makebox(0,0)[lb]{\smash{\SetFigFont{9}{10.8}{\familydefault}{\mddefault}{\updefault}{\color[rgb]{0,0,0}$c$}%
}}}
\put(6407,-1745){\makebox(0,0)[lb]{\smash{\SetFigFont{9}{10.8}{\familydefault}{\mddefault}{\updefault}{\color[rgb]{0,0,0}$c$}%
}}}
\put(6667,-1325){\makebox(0,0)[lb]{\smash{\SetFigFont{9}{10.8}{\familydefault}{\mddefault}{\updefault}{\color[rgb]{0,0,0}$c$}%
}}}
\put(7022,-985){\makebox(0,0)[lb]{\smash{\SetFigFont{9}{10.8}{\familydefault}{\mddefault}{\updefault}{\color[rgb]{0,0,0}$c^\ast$}%
}}}
\put(7467,-740){\makebox(0,0)[lb]{\smash{\SetFigFont{9}{10.8}{\familydefault}{\mddefault}{\updefault}{\color[rgb]{0,0,0}$c^\ast$}%
}}}
\put(6962,-4025){\makebox(0,0)[lb]{\smash{\SetFigFont{9}{10.8}{\familydefault}{\mddefault}{\updefault}{\color[rgb]{0,0,0}$c$}%
}}}
\put(6571,-3675){\makebox(0,0)[lb]{\smash{\SetFigFont{9}{10.8}{\familydefault}{\mddefault}{\updefault}{\color[rgb]{0,0,0}$c^\ast$}%
}}}
\put(6328,-3227){\makebox(0,0)[lb]{\smash{\SetFigFont{9}{10.8}{\familydefault}{\mddefault}{\updefault}{\color[rgb]{0,0,0}$c^\ast$}%
}}}
\put(6238,-2740){\makebox(0,0)[lb]{\smash{\SetFigFont{9}{10.8}{\familydefault}{\mddefault}{\updefault}{\color[rgb]{0,0,0}$c^\ast$}%
}}}
\end{picture}

%% file: fig4.pstex_t
\begin{picture}(0,0)%
\includegraphics{fig4.pstex}%
\end{picture}%
\setlength{\unitlength}{3947sp}%
\begingroup\makeatletter\ifx\SetFigFont\undefined%
\gdef\SetFigFont#1#2#3#4#5{%
  \reset@font\fontsize{#1}{#2pt}%
  \fontfamily{#3}\fontseries{#4}\fontshape{#5}%
  \selectfont}%
\fi\endgroup%
\begin{picture}(7282,2777)(901,-3180)
\put(901,-1711){\makebox(0,0)[lb]{\smash{\SetFigFont{10}{12.0}{\familydefault}{\mddefault}{\updefault}{\color[rgb]{0,0,0}$j=$}%
}}}
\put(901,-511){\makebox(0,0)[lb]{\smash{\SetFigFont{10}{12.0}{\familydefault}{\mddefault}{\updefault}{\color[rgb]{0,0,0}$k=$}%
}}}
\put(7351,-511){\makebox(0,0)[lb]{\smash{\SetFigFont{10}{12.0}{\familydefault}{\mddefault}{\updefault}{\color[rgb]{0,0,0}$4$}%
}}}
\put(5551,-511){\makebox(0,0)[lb]{\smash{\SetFigFont{10}{12.0}{\familydefault}{\mddefault}{\updefault}{\color[rgb]{0,0,0}$3$}%
}}}
\put(3751,-511){\makebox(0,0)[lb]{\smash{\SetFigFont{10}{12.0}{\familydefault}{\mddefault}{\updefault}{\color[rgb]{0,0,0}$2$}%
}}}
\put(1951,-511){\makebox(0,0)[lb]{\smash{\SetFigFont{10}{12.0}{\familydefault}{\mddefault}{\updefault}{\color[rgb]{0,0,0}$1$}%
}}}
\put(1198,-1572){\makebox(0,0)[lb]{\smash{\SetFigFont{12}{14.4}{\familydefault}{\mddefault}{\updefault}{\color[rgb]{0,0,0}$c$}%
}}}
\put(1498,-1572){\makebox(0,0)[lb]{\smash{\SetFigFont{12}{14.4}{\familydefault}{\mddefault}{\updefault}{\color[rgb]{0,0,0}$c$}%
}}}
\put(1798,-1572){\makebox(0,0)[lb]{\smash{\SetFigFont{12}{14.4}{\familydefault}{\mddefault}{\updefault}{\color[rgb]{0,0,0}$c$}%
}}}
\put(2098,-1572){\makebox(0,0)[lb]{\smash{\SetFigFont{12}{14.4}{\familydefault}{\mddefault}{\updefault}{\color[rgb]{0,0,0}$c^\ast$}%
}}}
\put(2398,-1572){\makebox(0,0)[lb]{\smash{\SetFigFont{12}{14.4}{\familydefault}{\mddefault}{\updefault}{\color[rgb]{0,0,0}$c^\ast$}%
}}}
\put(2698,-1572){\makebox(0,0)[lb]{\smash{\SetFigFont{12}{14.4}{\familydefault}{\mddefault}{\updefault}{\color[rgb]{0,0,0}$c^\ast$}%
}}}
\put(2998,-1572){\makebox(0,0)[lb]{\smash{\SetFigFont{12}{14.4}{\familydefault}{\mddefault}{\updefault}{\color[rgb]{0,0,0}$c$}%
}}}
\put(3298,-1572){\makebox(0,0)[lb]{\smash{\SetFigFont{12}{14.4}{\familydefault}{\mddefault}{\updefault}{\color[rgb]{0,0,0}$c$}%
}}}
\put(3598,-1572){\makebox(0,0)[lb]{\smash{\SetFigFont{12}{14.4}{\familydefault}{\mddefault}{\updefault}{\color[rgb]{0,0,0}$c$}%
}}}
\put(3898,-1572){\makebox(0,0)[lb]{\smash{\SetFigFont{12}{14.4}{\familydefault}{\mddefault}{\updefault}{\color[rgb]{0,0,0}$c^\ast$}%
}}}
\put(4198,-1572){\makebox(0,0)[lb]{\smash{\SetFigFont{12}{14.4}{\familydefault}{\mddefault}{\updefault}{\color[rgb]{0,0,0}$c^\ast$}%
}}}
\put(4498,-1572){\makebox(0,0)[lb]{\smash{\SetFigFont{12}{14.4}{\familydefault}{\mddefault}{\updefault}{\color[rgb]{0,0,0}$c^\ast$}%
}}}
\put(4798,-1572){\makebox(0,0)[lb]{\smash{\SetFigFont{12}{14.4}{\familydefault}{\mddefault}{\updefault}{\color[rgb]{0,0,0}$c$}%
}}}
\put(5098,-1572){\makebox(0,0)[lb]{\smash{\SetFigFont{12}{14.4}{\familydefault}{\mddefault}{\updefault}{\color[rgb]{0,0,0}$c$}%
}}}
\put(5398,-1572){\makebox(0,0)[lb]{\smash{\SetFigFont{12}{14.4}{\familydefault}{\mddefault}{\updefault}{\color[rgb]{0,0,0}$c$}%
}}}
\put(5698,-1572){\makebox(0,0)[lb]{\smash{\SetFigFont{12}{14.4}{\familydefault}{\mddefault}{\updefault}{\color[rgb]{0,0,0}$c^\ast$}%
}}}
\put(5998,-1572){\makebox(0,0)[lb]{\smash{\SetFigFont{12}{14.4}{\familydefault}{\mddefault}{\updefault}{\color[rgb]{0,0,0}$c^\ast$}%
}}}
\put(6298,-1572){\makebox(0,0)[lb]{\smash{\SetFigFont{12}{14.4}{\familydefault}{\mddefault}{\updefault}{\color[rgb]{0,0,0}$c^\ast$}%
}}}
\put(6598,-1572){\makebox(0,0)[lb]{\smash{\SetFigFont{12}{14.4}{\familydefault}{\mddefault}{\updefault}{\color[rgb]{0,0,0}$c$}%
}}}
\put(6898,-1572){\makebox(0,0)[lb]{\smash{\SetFigFont{12}{14.4}{\familydefault}{\mddefault}{\updefault}{\color[rgb]{0,0,0}$c$}%
}}}
\put(7198,-1572){\makebox(0,0)[lb]{\smash{\SetFigFont{12}{14.4}{\familydefault}{\mddefault}{\updefault}{\color[rgb]{0,0,0}$c$}%
}}}
\put(7498,-1572){\makebox(0,0)[lb]{\smash{\SetFigFont{12}{14.4}{\familydefault}{\mddefault}{\updefault}{\color[rgb]{0,0,0}$c^\ast$}%
}}}
\put(7798,-1572){\makebox(0,0)[lb]{\smash{\SetFigFont{12}{14.4}{\familydefault}{\mddefault}{\updefault}{\color[rgb]{0,0,0}$c^\ast$}%
}}}
\put(8098,-1572){\makebox(0,0)[lb]{\smash{\SetFigFont{12}{14.4}{\familydefault}{\mddefault}{\updefault}{\color[rgb]{0,0,0}$c^\ast$}%
}}}
\put(1501,-1711){\makebox(0,0)[lb]{\smash{\SetFigFont{8}{9.6}{\familydefault}{\mddefault}{\updefault}{\color[rgb]{0,0,0}2}%
}}}
\put(2401,-1711){\makebox(0,0)[lb]{\smash{\SetFigFont{8}{9.6}{\familydefault}{\mddefault}{\updefault}{\color[rgb]{0,0,0}2}%
}}}
\put(1201,-1711){\makebox(0,0)[lb]{\smash{\SetFigFont{8}{9.6}{\familydefault}{\mddefault}{\updefault}{\color[rgb]{0,0,0}3}%
}}}
\put(1801,-1711){\makebox(0,0)[lb]{\smash{\SetFigFont{8}{9.6}{\familydefault}{\mddefault}{\updefault}{\color[rgb]{0,0,0}1}%
}}}
\put(2101,-1711){\makebox(0,0)[lb]{\smash{\SetFigFont{8}{9.6}{\familydefault}{\mddefault}{\updefault}{\color[rgb]{0,0,0}1}%
}}}
\put(2701,-1711){\makebox(0,0)[lb]{\smash{\SetFigFont{8}{9.6}{\familydefault}{\mddefault}{\updefault}{\color[rgb]{0,0,0}3}%
}}}
\put(3301,-1711){\makebox(0,0)[lb]{\smash{\SetFigFont{8}{9.6}{\familydefault}{\mddefault}{\updefault}{\color[rgb]{0,0,0}2}%
}}}
\put(4201,-1711){\makebox(0,0)[lb]{\smash{\SetFigFont{8}{9.6}{\familydefault}{\mddefault}{\updefault}{\color[rgb]{0,0,0}2}%
}}}
\put(3001,-1711){\makebox(0,0)[lb]{\smash{\SetFigFont{8}{9.6}{\familydefault}{\mddefault}{\updefault}{\color[rgb]{0,0,0}3}%
}}}
\put(3601,-1711){\makebox(0,0)[lb]{\smash{\SetFigFont{8}{9.6}{\familydefault}{\mddefault}{\updefault}{\color[rgb]{0,0,0}1}%
}}}
\put(5101,-1711){\makebox(0,0)[lb]{\smash{\SetFigFont{8}{9.6}{\familydefault}{\mddefault}{\updefault}{\color[rgb]{0,0,0}2}%
}}}
\put(6001,-1711){\makebox(0,0)[lb]{\smash{\SetFigFont{8}{9.6}{\familydefault}{\mddefault}{\updefault}{\color[rgb]{0,0,0}2}%
}}}
\put(3901,-1711){\makebox(0,0)[lb]{\smash{\SetFigFont{8}{9.6}{\familydefault}{\mddefault}{\updefault}{\color[rgb]{0,0,0}1}%
}}}
\put(4501,-1711){\makebox(0,0)[lb]{\smash{\SetFigFont{8}{9.6}{\familydefault}{\mddefault}{\updefault}{\color[rgb]{0,0,0}3}%
}}}
\put(4801,-1711){\makebox(0,0)[lb]{\smash{\SetFigFont{8}{9.6}{\familydefault}{\mddefault}{\updefault}{\color[rgb]{0,0,0}3}%
}}}
\put(5401,-1711){\makebox(0,0)[lb]{\smash{\SetFigFont{8}{9.6}{\familydefault}{\mddefault}{\updefault}{\color[rgb]{0,0,0}1}%
}}}
\put(6901,-1711){\makebox(0,0)[lb]{\smash{\SetFigFont{8}{9.6}{\familydefault}{\mddefault}{\updefault}{\color[rgb]{0,0,0}2}%
}}}
\put(7801,-1711){\makebox(0,0)[lb]{\smash{\SetFigFont{8}{9.6}{\familydefault}{\mddefault}{\updefault}{\color[rgb]{0,0,0}2}%
}}}
\put(5701,-1711){\makebox(0,0)[lb]{\smash{\SetFigFont{8}{9.6}{\familydefault}{\mddefault}{\updefault}{\color[rgb]{0,0,0}1}%
}}}
\put(6301,-1711){\makebox(0,0)[lb]{\smash{\SetFigFont{8}{9.6}{\familydefault}{\mddefault}{\updefault}{\color[rgb]{0,0,0}3}%
}}}
\put(6601,-1711){\makebox(0,0)[lb]{\smash{\SetFigFont{8}{9.6}{\familydefault}{\mddefault}{\updefault}{\color[rgb]{0,0,0}3}%
}}}
\put(7201,-1711){\makebox(0,0)[lb]{\smash{\SetFigFont{8}{9.6}{\familydefault}{\mddefault}{\updefault}{\color[rgb]{0,0,0}1}%
}}}
\put(7501,-1711){\makebox(0,0)[lb]{\smash{\SetFigFont{8}{9.6}{\familydefault}{\mddefault}{\updefault}{\color[rgb]{0,0,0}1}%
}}}
\put(8101,-1711){\makebox(0,0)[lb]{\smash{\SetFigFont{8}{9.6}{\familydefault}{\mddefault}{\updefault}{\color[rgb]{0,0,0}3}%
}}}
\end{picture}

%% file: fig5.pstex_t
\begin{picture}(0,0)%
\includegraphics{fig5.pstex}%
\end{picture}%
\setlength{\unitlength}{3947sp}%
\begingroup\makeatletter\ifx\SetFigFont\undefined%
\gdef\SetFigFont#1#2#3#4#5{%
  \reset@font\fontsize{#1}{#2pt}%
  \fontfamily{#3}\fontseries{#4}\fontshape{#5}%
  \selectfont}%
\fi\endgroup%
\begin{picture}(6024,2724)(1189,-2023)
\put(5626, 14){\makebox(0,0)[lb]{\smash{\SetFigFont{12}{14.4}{\familydefault}{\mddefault}{\updefault}{\color[rgb]{0,0,0}$\Phi_2=$}%
}}}
\put(5626,-1486){\makebox(0,0)[lb]{\smash{\SetFigFont{12}{14.4}{\familydefault}{\mddefault}{\updefault}{\color[rgb]{0,0,0}$\Phi_2=$}%
}}}
\put(5026, 14){\makebox(0,0)[lb]{\smash{\SetFigFont{12}{14.4}{\familydefault}{\mddefault}{\updefault}{\color[rgb]{0,0,0}$\Longrightarrow$}%
}}}
\put(5026,-1486){\makebox(0,0)[lb]{\smash{\SetFigFont{12}{14.4}{\familydefault}{\mddefault}{\updefault}{\color[rgb]{0,0,0}$\Longrightarrow$}%
}}}
\put(5626,464){\makebox(0,0)[lb]{\smash{\SetFigFont{12}{14.4}{\familydefault}{\mddefault}{\updefault}{\color[rgb]{0,0,0}$\Phi_3=$}%
}}}
\put(5626,-436){\makebox(0,0)[lb]{\smash{\SetFigFont{12}{14.4}{\familydefault}{\mddefault}{\updefault}{\color[rgb]{0,0,0}$\Phi_1=$}%
}}}
\put(5626,-1036){\makebox(0,0)[lb]{\smash{\SetFigFont{12}{14.4}{\familydefault}{\mddefault}{\updefault}{\color[rgb]{0,0,0}$\Phi_3=$}%
}}}
\put(5626,-1936){\makebox(0,0)[lb]{\smash{\SetFigFont{12}{14.4}{\familydefault}{\mddefault}{\updefault}{\color[rgb]{0,0,0}$\Phi_1=$}%
}}}
\end{picture}

%% file: fig6.pstex_t
\begin{picture}(0,0)%
\includegraphics{fig6.pstex}%
\end{picture}%
\setlength{\unitlength}{3947sp}%
\begingroup\makeatletter\ifx\SetFigFont\undefined%
\gdef\SetFigFont#1#2#3#4#5{%
  \reset@font\fontsize{#1}{#2pt}%
  \fontfamily{#3}\fontseries{#4}\fontshape{#5}%
  \selectfont}%
\fi\endgroup%
\begin{picture}(6172,1354)(1801,-1808)
\put(1801,-961){\makebox(0,0)[lb]{\smash{\SetFigFont{12}{14.4}{\familydefault}{\mddefault}{\updefault}{\color[rgb]{0,0,0}$1$}%
}}}
\put(2101,-961){\makebox(0,0)[lb]{\smash{\SetFigFont{12}{14.4}{\familydefault}{\mddefault}{\updefault}{\color[rgb]{0,0,0}$2$}%
}}}
\put(2401,-961){\makebox(0,0)[lb]{\smash{\SetFigFont{12}{14.4}{\familydefault}{\mddefault}{\updefault}{\color[rgb]{0,0,0}$3$}%
}}}
\put(2701,-961){\makebox(0,0)[lb]{\smash{\SetFigFont{12}{14.4}{\familydefault}{\mddefault}{\updefault}{\color[rgb]{0,0,0}$4$}%
}}}
\put(3001,-961){\makebox(0,0)[lb]{\smash{\SetFigFont{12}{14.4}{\familydefault}{\mddefault}{\updefault}{\color[rgb]{0,0,0}$5$}%
}}}
\put(3301,-961){\makebox(0,0)[lb]{\smash{\SetFigFont{12}{14.4}{\familydefault}{\mddefault}{\updefault}{\color[rgb]{0,0,0}$6$}%
}}}
\put(3601,-961){\makebox(0,0)[lb]{\smash{\SetFigFont{12}{14.4}{\familydefault}{\mddefault}{\updefault}{\color[rgb]{0,0,0}$7$}%
}}}
\put(3901,-961){\makebox(0,0)[lb]{\smash{\SetFigFont{12}{14.4}{\familydefault}{\mddefault}{\updefault}{\color[rgb]{0,0,0}$8$}%
}}}
\put(5176,-586){\makebox(0,0)[lb]{\smash{\SetFigFont{12}{14.4}{\familydefault}{\mddefault}{\updefault}{\color[rgb]{0,0,0}$1$}%
}}}
\put(5551,-586){\makebox(0,0)[lb]{\smash{\SetFigFont{12}{14.4}{\familydefault}{\mddefault}{\updefault}{\color[rgb]{0,0,0}$2$}%
}}}
\put(5926,-586){\makebox(0,0)[lb]{\smash{\SetFigFont{12}{14.4}{\familydefault}{\mddefault}{\updefault}{\color[rgb]{0,0,0}$3$}%
}}}
\put(6301,-586){\makebox(0,0)[lb]{\smash{\SetFigFont{12}{14.4}{\familydefault}{\mddefault}{\updefault}{\color[rgb]{0,0,0}$4$}%
}}}
\put(6676,-586){\makebox(0,0)[lb]{\smash{\SetFigFont{12}{14.4}{\familydefault}{\mddefault}{\updefault}{\color[rgb]{0,0,0}$5$}%
}}}
\put(7051,-586){\makebox(0,0)[lb]{\smash{\SetFigFont{12}{14.4}{\familydefault}{\mddefault}{\updefault}{\color[rgb]{0,0,0}$6$}%
}}}
\put(7426,-586){\makebox(0,0)[lb]{\smash{\SetFigFont{12}{14.4}{\familydefault}{\mddefault}{\updefault}{\color[rgb]{0,0,0}$7$}%
}}}
\put(7801,-586){\makebox(0,0)[lb]{\smash{\SetFigFont{12}{14.4}{\familydefault}{\mddefault}{\updefault}{\color[rgb]{0,0,0}$8$}%
}}}
\put(4351,-1261){\makebox(0,0)[lb]{\smash{\SetFigFont{12}{14.4}{\familydefault}{\mddefault}{\updefault}{\color[rgb]{0,0,0}$\Longrightarrow$}%
}}}
\end{picture}

%% file: nested.pstex_t
\begin{picture}(0,0)%
\includegraphics{nested.pstex}%
\end{picture}%
\setlength{\unitlength}{3947sp}%
\begingroup\makeatletter\ifx\SetFigFont\undefined%
\gdef\SetFigFont#1#2#3#4#5{%
  \reset@font\fontsize{#1}{#2pt}%
  \fontfamily{#3}\fontseries{#4}\fontshape{#5}%
  \selectfont}%
\fi\endgroup%
\begin{picture}(2960,650)(1201,-526)
\put(1201,-286){\makebox(0,0)[lb]{\smash{\SetFigFont{12}{14.4}{\familydefault}{\mddefault}{\updefault}{\color[rgb]{0,0,0}$\varpi\,\; =$}%
}}}
\end{picture}

%% file: fig7.pstex_t
\begin{picture}(0,0)%
\includegraphics{fig7.pstex}%
\end{picture}%
\setlength{\unitlength}{3947sp}%
\begingroup\makeatletter\ifx\SetFigFont\undefined%
\gdef\SetFigFont#1#2#3#4#5{%
  \reset@font\fontsize{#1}{#2pt}%
  \fontfamily{#3}\fontseries{#4}\fontshape{#5}%
  \selectfont}%
\fi\endgroup%
\begin{picture}(3912,1869)(601,-1498)
\put(2401,239){\makebox(0,0)[lb]{\smash{\SetFigFont{12}{14.4}{\familydefault}{\mddefault}{\updefault}{\color[rgb]{0,0,0}$a$}%
}}}
\put(2701,239){\makebox(0,0)[lb]{\smash{\SetFigFont{12}{14.4}{\familydefault}{\mddefault}{\updefault}{\color[rgb]{0,0,0}$a$}%
}}}
\put(3001,239){\makebox(0,0)[lb]{\smash{\SetFigFont{12}{14.4}{\familydefault}{\mddefault}{\updefault}{\color[rgb]{0,0,0}$a^\ast$}%
}}}
\put(3301,239){\makebox(0,0)[lb]{\smash{\SetFigFont{12}{14.4}{\familydefault}{\mddefault}{\updefault}{\color[rgb]{0,0,0}$a^\ast$}%
}}}
\put(3601,-811){\makebox(0,0)[lb]{\smash{\SetFigFont{12}{14.4}{\familydefault}{\mddefault}{\updefault}{\color[rgb]{0,0,0}$a$}%
}}}
\put(3901,-811){\makebox(0,0)[lb]{\smash{\SetFigFont{12}{14.4}{\familydefault}{\mddefault}{\updefault}{\color[rgb]{0,0,0}$a$}%
}}}
\put(4201,-811){\makebox(0,0)[lb]{\smash{\SetFigFont{12}{14.4}{\familydefault}{\mddefault}{\updefault}{\color[rgb]{0,0,0}$a^\ast$}%
}}}
\put(4501,-811){\makebox(0,0)[lb]{\smash{\SetFigFont{12}{14.4}{\familydefault}{\mddefault}{\updefault}{\color[rgb]{0,0,0}$a^\ast$}%
}}}
\put(1201,-811){\makebox(0,0)[lb]{\smash{\SetFigFont{12}{14.4}{\familydefault}{\mddefault}{\updefault}{\color[rgb]{0,0,0}$a$}%
}}}
\put(1501,-811){\makebox(0,0)[lb]{\smash{\SetFigFont{12}{14.4}{\familydefault}{\mddefault}{\updefault}{\color[rgb]{0,0,0}$a$}%
}}}
\put(1801,-811){\makebox(0,0)[lb]{\smash{\SetFigFont{12}{14.4}{\familydefault}{\mddefault}{\updefault}{\color[rgb]{0,0,0}$a^\ast$}%
}}}
\put(2101,-811){\makebox(0,0)[lb]{\smash{\SetFigFont{12}{14.4}{\familydefault}{\mddefault}{\updefault}{\color[rgb]{0,0,0}$a^\ast$}%
}}}
\put(3601,239){\makebox(0,0)[lb]{\smash{\SetFigFont{12}{14.4}{\familydefault}{\mddefault}{\updefault}{\color[rgb]{0,0,0}$a$}%
}}}
\put(3901,239){\makebox(0,0)[lb]{\smash{\SetFigFont{12}{14.4}{\familydefault}{\mddefault}{\updefault}{\color[rgb]{0,0,0}$a$}%
}}}
\put(4201,239){\makebox(0,0)[lb]{\smash{\SetFigFont{12}{14.4}{\familydefault}{\mddefault}{\updefault}{\color[rgb]{0,0,0}$a^\ast$}%
}}}
\put(4501,239){\makebox(0,0)[lb]{\smash{\SetFigFont{12}{14.4}{\familydefault}{\mddefault}{\updefault}{\color[rgb]{0,0,0}$a^\ast$}%
}}}
\put(1201,239){\makebox(0,0)[lb]{\smash{\SetFigFont{12}{14.4}{\familydefault}{\mddefault}{\updefault}{\color[rgb]{0,0,0}$a$}%
}}}
\put(1501,239){\makebox(0,0)[lb]{\smash{\SetFigFont{12}{14.4}{\familydefault}{\mddefault}{\updefault}{\color[rgb]{0,0,0}$a$}%
}}}
\put(1801,239){\makebox(0,0)[lb]{\smash{\SetFigFont{12}{14.4}{\familydefault}{\mddefault}{\updefault}{\color[rgb]{0,0,0}$a^\ast$}%
}}}
\put(2101,239){\makebox(0,0)[lb]{\smash{\SetFigFont{12}{14.4}{\familydefault}{\mddefault}{\updefault}{\color[rgb]{0,0,0}$a^\ast$}%
}}}
\put(2401,-811){\makebox(0,0)[lb]{\smash{\SetFigFont{12}{14.4}{\familydefault}{\mddefault}{\updefault}{\color[rgb]{0,0,0}$a$}%
}}}
\put(2701,-811){\makebox(0,0)[lb]{\smash{\SetFigFont{12}{14.4}{\familydefault}{\mddefault}{\updefault}{\color[rgb]{0,0,0}$a$}%
}}}
\put(3001,-811){\makebox(0,0)[lb]{\smash{\SetFigFont{12}{14.4}{\familydefault}{\mddefault}{\updefault}{\color[rgb]{0,0,0}$a^\ast$}%
}}}
\put(3301,-811){\makebox(0,0)[lb]{\smash{\SetFigFont{12}{14.4}{\familydefault}{\mddefault}{\updefault}{\color[rgb]{0,0,0}$a^\ast$}%
}}}
\put(601,-61){\makebox(0,0)[lb]{\smash{\SetFigFont{12}{14.4}{\familydefault}{\mddefault}{\updefault}{\color[rgb]{0,0,0}$\pi$}%
}}}
\put(601,-1111){\makebox(0,0)[lb]{\smash{\SetFigFont{12}{14.4}{\familydefault}{\mddefault}{\updefault}{\color[rgb]{0,0,0}$\pi_r$}%
}}}
\end{picture}